\documentclass[11pt,reqno]{amsart}

\usepackage{setspace}
\usepackage[margin=1.1in]{geometry}
\usepackage{amsmath,amssymb,amsthm,mathtools}
\usepackage{graphicx}
\usepackage{booktabs}
\usepackage{tikz}
\usetikzlibrary{arrows.meta,positioning}
\usepackage{microtype}
\usepackage[colorlinks,linkcolor=blue!60!black,citecolor=blue!60!black,urlcolor=blue!60!black]{hyperref}
\usepackage[capitalise,noabbrev]{cleveref}

\theoremstyle{plain}
\newtheorem{theorem}{Theorem}[section]
\newtheorem{lemma}[theorem]{Lemma}

\theoremstyle{definition}
\newtheorem{definition}[theorem]{Definition}

\theoremstyle{remark}
\newtheorem{remark}[theorem]{Remark}

\numberwithin{equation}{section}

\title[Discrete Approximation of differential two-stage SLCP]{Discrete Approximation of differential two-stage stochastic linear complementarity problems}

\author{Xiaojun Chen}
\address{Department of Applied Mathematics,
The Hong Kong Polytechnic University, Kowloon, Hong Kong}
\email{xiaojun.chen@polyu.edu.hk}

\author{Jian Guo}
\address{Department of Applied Mathematics,
The Hong Kong Polytechnic University, Kowloon, Hong Kong}
\email{jians.guo@polyu.edu.hk}

\author{Seungjun Lee}
\address{School of Mathematics and Computing
(Computational Science and Engineering),
Yonsei University, Seoul 03722, Republic of Korea}
\email{sjlee0088@gmail.com}

\thanks{This work is supported by RGC grant JLFS/P-501/24
for the CAS AMSS--PolyU Joint Laboratory in Applied Mathematics
and Hong Kong Research Grants Council projects
PolyU15300123 and PolyU15300124.}

\subjclass[2020]{Primary 90C33;
Secondary 90C15, 90C31}
\keywords{Differential variational inequality, two-stage stochastic linear complementarity problem, time-dependent distribution, classical solution, progressive hedging method, Euler method}

\date{\today}

\usepackage{amsfonts}
\usepackage{graphicx}
\usepackage{epstopdf}
\usepackage{algorithmic}
\usepackage{algorithm}
\usepackage{booktabs}
\usepackage{enumerate}

\newcommand{\mb}{\mathbb}
\newcommand{\mc}{\mathcal}
\newcommand{\mf}{\mathbf}
\newcommand{\ip}[2]{\left\langle #1, #2 \right\rangle}
\newcommand{\lcp}{\mathrm{LCP}}
\newcommand{\sol}{\mathrm{SOL}}
\newcommand{\tv}{\mathrm{TV}}

\newtheorem{assumption}{Assumption}[section]

\newcommand{\nhe}{N,h,\eta}
\renewcommand{\theequation}{\thesection.\arabic{equation}}

\makeatletter

\def\@settitle{%
  \vspace{0.5em}
  {\centering
   \LARGE\bfseries
   \@title\par}
  \vspace{0.5em}
}

\def\@setauthors{%
  \begingroup
  \def\thanks{\protect\thankswarning}%
  \trivlist
  \centering\normalsize
  \@topsep10\p@\relax
  \item\relax
  \author@andify\authors
  \def\\{\protect\linebreak}%
  \authors
  \endtrivlist
  \endgroup
}

\makeatother

\begin{document}

\begin{abstract}
This paper studies differential two-stage stochastic linear
complementarity problems, in which an ordinary differential
equation (ODE) is coupled with a two-stage stochastic linear
complementarity problem (SLCP) under a time-dependent probability
distribution. The right-hand side of the ODE involves a solution
pair (here-and-now and wait-and-see decisions) of the two-stage
SLCP, while the two-stage SLCP involves two functions depending on
the time $t$ and the state $x(t)$ of the ODE. Under appropriate
conditions, we prove that the solution function of the two-stage
SLCP is Lipschitz continuous in $x(t)$ and locally
$\frac{1}{2}$-H\"older continuous in $t$, and establish the existence and
uniqueness of a classical solution to the differential two-stage
SLCP. Moreover, we develop a discrete approximation scheme for the
differential two-stage SLCP. We first use a finite-scenario
partition with  mesh size $\delta_N$ to approximate the distribution, then apply the
Euler method for time discretization with step size $h$, and finally solve the
resulting finite-dimensional problems using an inexact
progressive hedging method with tolerance $\eta$.
We prove that the error between a computed solution by the discrete approximation scheme and the classical solution of the original
problem is bounded by $C(\delta_N+h+\eta)$, where $C>0$ is a constant independent of $N$, $h$ and $\eta$. Finally, we apply our
method to find an equilibrium point of a differential two-stage
stochastic two-player zero-sum game. Numerical experiments
illustrate the convergence rates and the efficiency of the discrete
approximation scheme.
\end{abstract}

\maketitle

\section{Introduction}\label{sec:introduction}
We consider the following differential two-stage stochastic linear complementarity problem (SLCP)
\begin{subequations}\label{DSVI}
\begin{align}
\dot{x}(t)
&= \mb{E}_{P_t}\big[\Phi\big(t,x(t),y(t),u(t,\xi),\xi\big)\big],
\qquad t\in[0,T],
\label{DSVI-1}\\
x(0)
&= \mf x_0,
\label{DSVI-2}\\
0 \leq y(t)
&\perp M_{11}y(t)
+ \mb{E}_{P_t}[M_{12}(\xi)u(t,\xi)]
+ f(t,x(t)) \geq 0,
\label{DSVI-3}\\
0 \leq u(t,\xi)
&\perp M_{21}(\xi)y(t)
+ M_{22}(\xi)u(t,\xi)
+ g(t,x(t),\xi) \geq 0,
\quad \text{$P_t$-a.e. } \xi\in\Xi.
\label{DSVI-4}
\end{align}
\end{subequations}
At each time $t\in[0,T]$, uncertainty is modeled by a random vector $\xi:\Omega\to\mb{R}^\ell$ on a probability space $(\Omega,\mc{F},\mb{P}_t)$ with time-dependent distribution $P_t:=\mb{P}_t\circ\xi^{-1}$ supported in a fixed Borel set $\Xi\subset\mb{R}^\ell$. There is a $\sigma$-finite Borel measure $\mu$ on $\Xi$ with $P_t\ll\mu$ for all $t$. The data are
$\Phi:\mb{R}_+\times\mb{R}^n\times\mb{R}^m\times\mb{R}^k\times\Xi\to\mb{R}^n$,
$M_{11}\in\mb{R}^{m\times m}$, $M_{12}:\Xi\to\mb{R}^{m\times k}$, $M_{21}:\Xi\to\mb{R}^{k\times m}$, $M_{22}:\Xi\to\mb{R}^{k\times k}$,
$f:\mb{R}_+\times\mb{R}^n\to\mb{R}^m$, $g:\mb{R}_+\times\mb{R}^n\times\Xi\to\mb{R}^k$. Here $x(t)$ is the state of the ordinary differential equation (ODE) in
\eqref{DSVI-1}--\eqref{DSVI-2}, and
$(y(t),u(t,\cdot))$ is a solution pair (here-and-now
and wait-and-see decisions) of the
two-stage SLCP in \eqref{DSVI-3}--\eqref{DSVI-4}
at a fixed $(t,x(t))$.

The problem defined by \eqref{DSVI-1}--\eqref{DSVI-3} with
$M_{12}(\xi)\equiv0$ and $u(t,\xi)\equiv0$ reduces to the
differential LCP \cite{PangStewart2008}, and the problem defined by
\eqref{DSVI-1}, \eqref{DSVI-2} and \eqref{DSVI-4} with
$M_{21}(\xi)\equiv0$, $y(t)\equiv0$ and $P_t\equiv P$ reduces to the
differential SLCP \cite{ChenShen2022}.
The ODE in most existing studies for differential variational
inequalities (DVIs) involves either a here-and-now solution or a
wait-and-see solution of the VI
\cite{ChenICM2026,chen2025differential,ChenShen2022,PangStewart2008}.
In contrast, the function $\Phi$ in \eqref{DSVI-1} involves a solution pair
(here-and-now and wait-and-see
decisions) of the two-stage SLCP.

The differential two-stage SLCP \eqref{DSVI} is motivated
by systems in which advance decisions and later adjustments
jointly affect the state.
For example, in energy storage scheduling, day-ahead
energy commitments and hour-ahead charging or discharging
adjustments both enter the battery's energy balance~\cite{AkhavanHejaziMohsenianRad2014}.
This motivates including both the here-and-now decision
$y(t)$ and the wait-and-see decision $u(t,\xi)$
in the ODE \eqref{DSVI-1}.
These decisions solve the two-stage SLCP
\eqref{DSVI-3}--\eqref{DSVI-4} at the current state $x(t)$,
with the two-stage structure studied in
\cite{chen2019discrete,ZhangSunXu2019}.

A second motivation is that the probability distribution
used to make decisions may change over time.
In energy storage scheduling, wind forecasts are updated
as new information becomes available
\cite{GhadimiPowell2024}.
In online inventory routing, new demand data are used
to update scenario distributions and plan deliveries
\cite{JiaSchrotenboerChen2025}.
Such updates not only affect the decisions themselves,
but also the expected state response.
In the history-dependent model of
\cite{chen2026stability}, the probability law enters
the expected state dynamics, while the lower-level
response is determined at each fixed state and sample path.
In \eqref{DSVI}, $P_t$ enters both the ODE \eqref{DSVI-1}
and, through the expected recourse term
$\mb E_{P_t}[M_{12}(\xi)u(t,\xi)]$ in \eqref{DSVI-3},
the two-stage SLCP.
Thus, changes in $P_t$ can alter
the here-and-now decision and its associated recourse
even at a fixed state.

We review some related work on two-stage SLCPs and DVIs.

Two-stage SLCPs are a class of two-stage stochastic
variational inequalities (SVIs) and have been extensively
studied
\cite{ChenPongWets2017,ChenShapiroSun2019,
chen2019discrete,RW2017,SunChen2021}.
They arise as optimality conditions of two-stage
stochastic linear programs with recourse
\cite{BirgeLouveaux2011,RW1976recourse,SDR2009,Wets1974}
and as equilibrium conditions of two-stage stochastic games
\cite{ChenPongWets2017,RW2017,SunChen2021}.
Two-stage  SVIs
arising from stochastic programming were studied
in \cite{LiZhang2020}.
A formulation of two-stage and multistage SVIs was developed in \cite{RW2017}, and the
progressive hedging method \cite{rockafellar1991scenarios} was
extended to monotone SVIs and
complementarity problems with finitely many scenarios in
\cite{rockafellar2019solving}.
Related studies also consider stochastic approximation methods
\cite{ChenLiuYangZhang2022,YangZhangWangLin2021}, nonmonotone SVIs \cite{CuiZhang2024}, and robust formulations
\cite{HoriYamashita2022,SunShapiroChen2023}.
See \cite{SunChen2021} for a survey.
Applications include Cournot--Nash competition
\cite{JiangShiWangChen2019}, electricity markets
\cite{DeGiuliMilasiOggioniScopelliti2026}, and the storage and
distribution of medical supplies \cite{LiZhangDingLv2022}.
More recently, \cite{SunChen2025} showed that the KKT
conditions of a class of two-stage stochastic minimax
problems can be expressed as a two-stage SVI.

DVIs couple ODEs with parametric variational inequalities \cite{PangStewart2008}. They include differential complementarity systems and provide a framework for studying dynamics subject to equilibrium constraints. For linear complementarity systems, Zeno behavior has been studied in \cite{CamlibelPangShen2006,Shen2014,ShenPang2005,ShenPang2007}. Numerical methods include time-stepping schemes based on least-norm selections \cite{Han2009}, Newton iterations in implicit time stepping \cite{chen2013newton}, and Tikhonov regularization \cite{chen2013convergence}.
Stochastic extensions include dynamic stochastic variational inequalities (DSVIs) driven by expected scenario-wise
VI solutions, with results on existence, uniqueness, uniform sample
average approximation (SAA)
convergence with an exponential rate, and time-stepping convergence
\cite{ChenShen2022}. History-dependent stochastic responses were further
studied, with results on well-posedness, uniform SAA convergence,
and Wasserstein-based transfer stability \cite{chen2026stability}.

Replacing the VIs in DVIs by stochastic optimization provides another class of problems related to \eqref{DSVI-1}--\eqref{DSVI-4}. In \cite{luo2025dynamic}, the function in the ODE involves a here-and-now solution of a stochastic nonsmooth convex optimization problem, and in \cite{chen2025differential}, the function in the ODE involves wait-and-see solutions of scenario-wise convex optimization problems and permits time-dependent probability distributions.

{\bf Contributions.}
The paper has three main contributions.
\begin{enumerate}[(i)]
\item
We establish Lipschitz continuity in the state and local $\frac{1}{2}$-H\"older continuity in time for the solution of the two-stage SLCP (Lemma~\ref{lem:est}). Based on these results, we prove the existence and uniqueness of a classical solution to the differential two-stage SLCP, with a continuously differentiable state, a continuous here-and-now decision, and a jointly measurable wait-and-see decision (Theorem~\ref{thm:classical}).
We also provide a counterexample showing that separate moment bounds do not ensure the existence of classical solutions (Remark~\ref{rem:joint}).

\item We develop a discrete approximation scheme for the differential two-stage SLCP. We first construct a finite-scenario approximation of the
differential two-stage SLCP by partitioning the uncertainty set
and establish an $O(\delta_N)$ convergence rate, where
$\delta_N$ is the mesh size of the partition
(Theorem~\ref{thm:semi-discrete}).
We then apply the forward Euler method with time step $h$
to the finite-scenario problem and solve the resulting
finite-dimensional problems using an inexact progressive
hedging method with stopping tolerance $\eta$.
We prove that the error between a computed solution by the discrete approximation scheme and the classical solution of~\eqref{DSVI} is bounded by
$C(\delta_N+h+\eta)$, where $C$ is independent of
$N$, $h$ and $\eta$ (Theorem~\ref{thm:euler}).

\item
We formulate a differential two-stage SLCP from the KKT conditions
of a quadratic two-stage stochastic two-player
zero-sum game \cite{SunChen2025}.
Its equilibria correspond to the KKT points
of the associated minimax problem.
We consider examples with nonunique and unique complementarity
solutions and investigate the role of regularization.
Numerical experiments support the first-order convergence rates
under strong monotonicity and show similar behavior when
the solution is unique but strong monotonicity fails
(Section~\ref{sec:minimax-experiments}).

\end{enumerate}

{\bf Organization of the paper.}
The rest of this section fixes the notation, the solution concept and the standing assumptions. Section~\ref{sec:classical} establishes the existence and uniqueness of classical solutions. Section~\ref{sec:discretization} introduces the discretization and proves the error estimates. Section~\ref{sec:minimax-experiments} contains the application to two-stage stochastic minimax problems and the numerical experiments, and Section~\ref{sec:conclusion} concludes the paper.

\subsection{Notation}\label{subsec:notation}

Boldface symbols $\mf y,\mf u,\mf z$ denote elements
of the spaces introduced below (vectors or functions of $\xi$),
whereas nonbold symbols are used for state and decision
trajectories, with $x(t)$, $y(t)$, $u(t,\xi)$, and $z(t,\xi)$
denoting their values at the specified arguments.
For $N\in\mb N$, $[N]:=\{1,\dots,N\}$. For $\mf v\in\mb{R}^p$, $\|\mf v\|$ is the Euclidean norm. $\mf 0_p$ and $\mf 1_p$ are the all-zeros and the all-ones vectors of length $p$, and $I_p$ is the $p\times p$ identity matrix. Inequalities between vectors, such as $\mf v \geq \mf 0_p$, as well as $\max\{\cdot,\cdot\}$ and $\min\{\cdot,\cdot\}$ of vectors are understood component-wise. For a matrix, $\|\cdot\|$ is the spectral norm. A square matrix $A$ is positive definite if $\mf v^\top A\mf v>0$ for all $\mf v\ne\mf 0$, and for symmetric matrices $A,B$, $A\succeq B$ means that $A-B$ is positive semidefinite.

For each probability measure $P_t$ on $\Xi$, we define
the space
\[
\resizebox{1\linewidth}{!}{$\displaystyle
\mc H_{P_t}:=\mb R^m\times L^2_{P_t}(\mb R^k)
=\Big\{\mf z=(\mf y,\mf u):\ \mf y\in\mb R^m,\
\mf u:\Xi\to\mb R^k\ \text{measurable},\
\mb E_{P_t}\|\mf u(\xi)\|^2<\infty\Big\}
$}
\]
equipped with the norm
$
\|\mf z\|_{\mc H_{P_t}}^2
:=\|\mf y\|^2+\mb E_{P_t}\|\mf u(\xi)\|^2,
$
where functions in $L^2_{P_t}(\mb R^k)$
are identified if they agree $P_t$-a.s.
For $\mf z=(\mf y,\mf u)\in\mc H_{P_t}$, we write
$\mf z(\xi):=(\mf y,\mf u(\xi))$ and use $\Phi(t,\mf x,\mf z(\xi),\xi)$
as shorthand for $\Phi(t,\mf x,\mf y,\mf u(\xi),\xi)$.
When $\mf z$ depends on parameters, as in
$\mf z(t,\mf x)=(\mf y(t,\mf x),\mf u(t,\mf x))$, we write
$\mf u(t,\mf x,\xi)$ and $\mf z(t,\mf x,\xi)$ for the values at $\xi$. For trajectories, we write $z(t,\xi):=(y(t),u(t,\xi))$.

For \eqref{DSVI-3}--\eqref{DSVI-4}, we define the block matrix
$
M(\xi):=
\begin{pmatrix}
M_{11}&M_{12}(\xi)\\
M_{21}(\xi)&M_{22}(\xi)
\end{pmatrix}
\in\mb R^{(m+k)\times(m+k)}.
$

\subsection{Two-stage SLCP formulation and assumptions}
\label{subsec:lcp}

For $s\in[0,T]$, $\mf f\in\mb R^m$ and measurable $\mf g:\Xi\to\mb R^k$, writing $\mf q=(\mf f,\mf g)$,
the two-stage SLCP with data $(s,\mf q)$ is to find
$\mf z=(\mf y,\mf u)\in\mc H_{P_s}$ such that
\begin{equation}\label{eq:tslcp}
\left\{
\begin{aligned}
&0\le\mf y\perp M_{11}\mf y+\mb E_{P_s}[M_{12}(\xi)\mf u(\xi)]+\mf f\ge0,\\
&0\le\mf u(\xi)\perp M_{21}(\xi)\mf y+M_{22}(\xi)\mf u(\xi)+\mf g(\xi)\ge0
\quad\text{for $P_s$-a.e.\ }\xi\in\Xi.
\end{aligned}
\right.
\end{equation}
We denote its solution set by $\sol_{P_s}(\mf q)\subset\mc H_{P_s}$.
For each $t\in[0,T]$ and $\mf x\in\mb R^n$, let
$q(t,\mf x):=(f(t,\mf x),g(t,\mf x,\cdot))$ and
$\mc S(t,\mf x):=\sol_{P_t}(q(t,\mf x))$.
Thus, \eqref{DSVI-3}--\eqref{DSVI-4}
can be written as $z(t)\in\mc S(t,x(t))$.
Under the assumptions below, $\mc S(t,\mf x)$ consists of a single element
(Lemma~\ref{lem:wellposed}), which we denote by
$\mf z(t,\mf x)=(\mf y(t,\mf x),\mf u(t,\mf x))$: the two-stage SLCP
solution at time $t$ and state $\mf x$, regarded as a function of
$(t,\mf x)$, with second-stage values $\mf u(t,\mf x,\xi)$.
The decision $z(t)$ of a classical solution is then
$z(t)=\mf z(t,x(t))$ along the state trajectory.

We now define a classical solution of the differential two-stage SLCP \eqref{DSVI}.
\begin{definition}[Classical solution]\label{def:classical}
Let $T>0$. A pair $(x,z)$
with $x:[0,T]\to\mb R^n$ and
$z:t\mapsto z(t)=(y(t),u(t,\cdot))\in\mc H_{P_t}$
is a classical solution of \eqref{DSVI}
on $[0,T]$ if $x$ is continuously differentiable,
$y$ is continuous, $u$ admits a jointly measurable
representative, and
$\dot x(t)=\mb E_{P_t}[\Phi(t,x(t),z(t,\xi),\xi)]$
and $z(t)\in\mc S(t,x(t))$ for every $t\in[0,T]$,
with $x(0)=\mf x_0$.
\end{definition}

For a fixed $T>0$, we introduce the following assumptions
for the well-posedness analysis in Section~\ref{sec:classical}.

\begin{assumption}[Condition on $\Phi$]\label{ass:Phi}
$\Phi(t,\mf x,\mf z,\cdot)$ is measurable,
and there exist a constant $\bar\kappa_\Phi\ge0$
and a measurable function
$\kappa_\Phi:\Xi\to[0,\infty)$ with
$\sup_{t\in[0,T]}
\!\mb E_{P_t}[\kappa_\Phi(\xi)^2]$\\$
\le\bar\kappa_\Phi^2$
such that, for $\mu$-a.e.\ $\xi\in\Xi$,
\begin{align*}
\|\Phi(t,\mf x,\mf z,\xi)\|
&\le\kappa_\Phi(\xi)
\big(1+t+\|\mf x\|+\|\mf z\|\big),\\
\|\Phi(t,\mf x,\mf z,\xi)
-\Phi(s,\mf x',\mf z',\xi)\|
&\le\kappa_\Phi(\xi)
\big(|t-s|+\|\mf x-\mf x'\|
+\|\mf z-\mf z'\|\big),
\end{align*}
for all $t,s\in[0,T]$,
$\mf x,\mf x'\in\mb R^n$ and
$\mf z,\mf z'\in\mb R^{m+k}$.
\end{assumption}

\begin{assumption}[Condition on $M(\xi)$]
\label{ass:matrix}
$M(\cdot)$ is measurable, and there exists a measurable function
$\kappa:\Xi\to(0,\infty)$ such that, for $\mu$-a.e.\ $\xi\in\Xi$,
$\mf z^\top M(\xi)\mf z\ge\kappa(\xi)\|\mf z\|^2$
for all $\mf z\in\mb R^{m+k}$.
Moreover, there exist constants $\bar M,\underline\kappa>0$ such that
\[
\sup_{t\in[0,T]}\mb E_{P_t}\big[\|M(\xi)\|^2\big]\le\bar M^2,
\qquad
\inf_{t\in[0,T]}\mb E_{P_t}[\kappa(\xi)]\ge\underline\kappa.
\]
\end{assumption}
\begin{assumption}[Condition on $q$]\label{ass:q}
$q(t,\mf x,\cdot)$ is measurable,
and there exist a constant $\bar\kappa_q\ge0$
and a measurable function
$\kappa_q:\Xi\to[0,\infty)$ with
$\sup_{t\in[0,T]}
\mb E_{P_t}[\kappa_q(\xi)^2]
\le\bar\kappa_q^2$
such that, for $\mu$-a.e.\ $\xi\in\Xi$ and for all $t,s\in[0,T]$,
$\mf x,\mf x'\in\mb R^n$, we have
\begin{align*}
\|q(t,\mf x,\xi)\|
&\le\kappa_q(\xi)(1+t+\|\mf x\|),\\
\|q(t,\mf x,\xi)-q(s,\mf x',\xi)\|
&\le\kappa_q(\xi)
\big(|t-s|+\|\mf x-\mf x'\|\big).
\end{align*}

\end{assumption}

\begin{assumption}[Condition on $P_t$]\label{ass:kernel}
There exists a constant $L_P>0$ such that $t\mapsto P_t$ is Lipschitz continuous:
$
\|P_t-P_s\|_{\tv}\le L_P|t-s|,
$
for all $t,s\in[0,T]$,
where
$\|P_t-P_s\|_{\tv}:=\int_\Xi|p_t-p_s|\,d\mu$
is the total variation norm with $p_t:=dP_t/d\mu$.
\end{assumption}

\begin{assumption}[Joint moment conditions]
\label{ass:joint}
Let $\theta(\xi):=(\|M(\xi)\|+\kappa_q(\xi))/\kappa(\xi)$. There are constants
$\bar\theta,\bar\kappa_{\Phi,q}\ge0$ such that
\[
\sup_{t\in[0,T]}\mb E_{P_t}\big[(1+\|M(\xi)\|)^2\theta(\xi)^2\big]\le\bar\theta^2,
\qquad
\sup_{t\in[0,T]}\mb E_{P_t}\big[\kappa_\Phi(\xi)^2\theta(\xi)^2\big]\le\bar\kappa_{\Phi,q}^2.
\]
\end{assumption}

The second-stage solution is controlled scenario-wise by $\theta(\xi)$
(Lemma~\ref{lem:wellposed}), and the first condition in
Assumption~\ref{ass:joint} ensures that it belongs to $L^2_{P_t}(\mb R^k)$.
The factor $(1+\|M(\xi)\|)^2$ there and the second condition are needed
to compare expectations under different measures
(Lemmas~\ref{lem:est}(iii) and \ref{lem:F}(iii)). Both conditions hold if
$\|M\|$, $1/\kappa$ and either $\kappa_\Phi$ or $\kappa_q$ are bounded.
Remark~\ref{rem:joint} shows that they cannot be dropped or replaced by
separate second-moment bounds.

\begin{remark}[Null set]\label{rem:null}
The pointwise conditions in Assumptions~\ref{ass:Phi}, \ref{ass:matrix}
and \ref{ass:q} hold outside a common $\mu$-null set. Since $P_t\ll\mu$,
this set is $P_t$-null for every $t\in[0,T]$, so redefining the data on it
affects neither the expectations in \eqref{DSVI-1} and \eqref{DSVI-3} nor
the $P_t$-a.e.\ condition in \eqref{DSVI-4}, and hence not the notion of a
classical solution. Redefining the data on this set, we assume that these
pointwise conditions hold for every $\xi\in\Xi$.
\end{remark}

\section{Well-posedness of the differential two-stage SLCP}\label{sec:classical}

In this section,
we first prove the well-posedness of the two-stage SLCP
\eqref{eq:tslcp} and then derive the growth and continuity estimates
needed for the dynamic system with time-dependent
probability measures.
Throughout this section,
Assumptions~\ref{ass:Phi}--\ref{ass:joint} are in force.
Set $\rho(t,\mf x):=1+t+\|\mf x\|$.
\subsection{Basic estimates and the fixed-parameter SLCP}

For $\mf z=(\mf y,\mf u)\in\mc H_{P_t}$, since $P_t$ is a probability
measure,
\begin{equation}\label{eq:z}
\|\mf z\|_{\mc H_{P_t}}^2
=\|\mf y\|^2+\mb E_{P_t}\|\mf u(\xi)\|^2
=\mb E_{P_t}\|\mf z(\xi)\|^2.
\end{equation}

Using Proposition~2.1 in \cite{chen2019discrete}, we give the following lemma for the two-stage SLCP in~\eqref{eq:tslcp}.

\begin{lemma}[Well-posedness of the two-stage SLCP]
\label{lem:wellposed}
Let $s\in[0,T]$ and $\mf q=(\mf f,\mf g)$ be as in \eqref{eq:tslcp} with
$\mb E_{P_s}[\|\mf g(\xi)\|^2/\kappa(\xi)^2]<\infty$.
\begin{enumerate}[(i)]
\item For every $\mf y\in\mb R^m$ and $\xi\in\Xi$,
the second-stage LCP
\begin{equation}\label{eq:lcp-xi}
0\le\mf u\perp M_{22}(\xi)\mf u+M_{21}(\xi)\mf y+\mf g(\xi)\ge0
\end{equation}
has a unique solution $\mf u_{\mf y,\mf g}(\xi)$, measurable in $\xi$, and
for any other pair $(\mf y',\mf g')$,
\begin{equation}\label{eq:u-lip-y}
\begin{aligned}
\|\mf u_{\mf y,\mf g}(\xi)\|
&\le\frac{\|M_{21}(\xi)\|\|\mf y\|+\|\mf g(\xi)\|}{\kappa(\xi)},\\
\|\mf u_{\mf y,\mf g}(\xi)-\mf u_{\mf y',\mf g'}(\xi)\|
&\le\frac{\|M_{21}(\xi)\|\|\mf y-\mf y'\|
+\|\mf g(\xi)-\mf g'(\xi)\|}{\kappa(\xi)}.
\end{aligned}
\end{equation}
\item $\sol_{P_s}(\mf q)=\{(\mf y,\mf u_{\mf y,\mf g})\}$, where $\mf y$ is
the unique solution of the nonlinear complementarity
problem (NCP)
\begin{equation}\label{eq:reduced}
0\le\mf y\perp M_{11}\mf y+\mb E_{P_s}[M_{12}(\xi)\mf u_{\mf y,\mf g}(\xi)]
+\mf f\ge0,
\end{equation}
whose map is Lipschitz continuous and strongly monotone on $\mb R^m$ with
modulus $\mb E_{P_s}[\kappa(\xi)]\ge\underline\kappa$.
\item If $(\mf y',\mf u')\in\sol_{P_s}(\mf q')$ for another
$\mf q'=(\mf f',\mf g')$ as above, then
\begin{equation}\label{eq:perturb}
\|\mf y-\mf y'\|\le\frac1{\underline\kappa}\Big(\|\mf f-\mf f'\|
+\mb E_{P_s}\Big[\|M_{12}(\xi)\|
\frac{\|\mf g(\xi)-\mf g'(\xi)\|}{\kappa(\xi)}\Big]\Big).
\end{equation}
\end{enumerate}
\end{lemma}

\begin{proof}
(i) By Assumption~\ref{ass:matrix}, strong monotonicity
guarantees a unique solution to \eqref{eq:lcp-xi}. It is measurable in $\xi$, since
it coincides on finitely many measurable sets, one for each active set
$J\subseteq[k]$, with the measurable solution of the corresponding linear
system. For $\mf u:=\mf u_{\mf y,\mf g}(\xi)$ and
$\mf u':=\mf u_{\mf y',\mf g'}(\xi)$, adding the variational inequalities
$(\mf u'-\mf u)^\top(M_{22}(\xi)\mf u+M_{21}(\xi)\mf y+\mf g(\xi))\ge0$
and
$(\mf u-\mf u')^\top(M_{22}(\xi)\mf u'+M_{21}(\xi)\mf y'+\mf g'(\xi))\ge0$
gives
$\kappa(\xi)\|\mf u-\mf u'\|^2\le\|\mf u-\mf u'\|
\big(\|M_{21}(\xi)\|\|\mf y-\mf y'\|+\|\mf g(\xi)-\mf g'(\xi)\|\big)$,
which is the second bound in \eqref{eq:u-lip-y}. The first follows with
$\mf y'=\mf 0_m$ and $\mf g'=\max\{\mf g,\mf 0_k\}$, for which $\mf u'=\mf 0_k$.

\smallskip\noindent
(ii) Let
$G(\mf y):=M_{11}\mf y+\mb E_{P_s}[M_{12}(\xi)\mf u_{\mf y,\mf g}(\xi)]+\mf f$,
so that \eqref{eq:reduced} reads $0\le\mf y\perp G(\mf y)\ge0$. By
\eqref{eq:u-lip-y}, Assumptions~\ref{ass:matrix} and \ref{ass:joint} and the
Cauchy--Schwarz inequality, $\mf u_{\mf y,\mf g}\in L^2_{P_s}(\mb R^k)$ and
$G$ is well defined and Lipschitz continuous. Let $\mf d:=\mf y-\mf y'$ and
$\mf w:=\mf u_{\mf y,\mf g}(\xi)-\mf u_{\mf y',\mf g}(\xi)$. As in (i),
$\mf w^\top(M_{22}(\xi)\mf w+M_{21}(\xi)\mf d)\le0$, so Assumption~\ref{ass:matrix}
applied to $(\mf d,\mf w)$ yields, suppressing the dependence on $\xi$,
\[
\mf d^\top M_{11}\mf d+\mf d^\top M_{12}\mf w
\ge\kappa\big(\|\mf d\|^2+\|\mf w\|^2\big)
-\mf w^\top(M_{22}\mf w+M_{21}\mf d)\ge\kappa\|\mf d\|^2.
\]
Taking the expectation under $P_s$ gives
$\mf d^\top(G(\mf y)-G(\mf y'))\ge\mb E_{P_s}[\kappa]\|\mf d\|^2$.
Hence \eqref{eq:reduced} has a unique solution $\mf y^*$
\cite{kinderlehrer2000introduction}, and $(\mf y^*,\mf u_{\mf y^*,\mf g})$
satisfies \eqref{eq:tslcp}. Conversely, if $(\mf y,\mf u)\in\sol_{P_s}(\mf q)$,
then $\mf u=\mf u_{\mf y,\mf g}$ $P_s$-a.e.\ by (i), so $\mf y$ solves
\eqref{eq:reduced} and $\mf y=\mf y^*$.

\smallskip\noindent
(iii) Let $G'$ be the map in \eqref{eq:reduced} for $\mf q'$. Adding the
variational inequalities of the NCPs for $G$ and $G'$, and using the strong
monotonicity of $G$ and the second bound in \eqref{eq:u-lip-y} with
$\mf y=\mf y'$, we obtain
\[
\underline\kappa\|\mf y-\mf y'\|
\le\|G'(\mf y')-G(\mf y')\|
\le\|\mf f-\mf f'\|
+\mb E_{P_s}\Big[\|M_{12}(\xi)\|
\frac{\|\mf g(\xi)-\mf g'(\xi)\|}{\kappa(\xi)}\Big].
\]
\end{proof}

We now apply Lemma~\ref{lem:wellposed} to the data $q(t,\mf x)$ and
bound the two terms on the right-hand side of \eqref{eq:perturb}.
Since $\|f(t,\mf x)\|\le\|q(t,\mf x,\xi)\|\le\kappa_q(\xi)\rho(t,\mf x)$
for every $\xi\in\Xi$, and similarly for differences,
Assumption~\ref{ass:q} gives, for all $t,s\in[0,T]$ and
$\mf x,\mf x'\in\mb R^n$,
\begin{equation}\label{eq:q-bounds}
\|f(t,\mf x)\|\le\bar\kappa_q\rho(t,\mf x),\qquad
\|f(t,\mf x)-f(s,\mf x')\|\le\bar\kappa_q\big(|t-s|+\|\mf x-\mf x'\|\big).
\end{equation}
Moreover,
$\|g(t,\mf x,\xi)-g(s,\mf x',\xi)\|/\kappa(\xi)
\le\theta(\xi)(|t-s|+\|\mf x-\mf x'\|)$, so that
Assumption~\ref{ass:joint} gives, for every $r\in[0,T]$,
\begin{equation}\label{eq:g-ratio}
\mb E_{P_r}\Big[\big(1+\|M(\xi)\|\big)^2
\frac{\|g(t,\mf x,\xi)-g(s,\mf x',\xi)\|^2}{\kappa(\xi)^2}\Big]^{1/2}
\le\bar\theta\big(|t-s|+\|\mf x-\mf x'\|\big),
\end{equation}
and the same bound holds with $g(s,\mf x',\xi)$ replaced by $\mf 0_k$
and the right-hand side by $\bar\theta\rho(t,\mf x)$.
In particular, the integrability condition of
Lemma~\ref{lem:wellposed} holds for $q(t,\mf x)$, and
$\mc S(t,\mf x)$ contains exactly one element, denoted by
$\mf z(t,\mf x)=(\mf y(t,\mf x),\mf u(t,\mf x))\in\mc H_{P_t}$.
We use the representative
$\mf u(t,\mf x,\xi):=\mf u_{\mf y(t,\mf x),g(t,\mf x,\cdot)}(\xi)$,
defined for every $\xi\in\Xi$, and similarly for the element of
$\sol_{P_s}(\mf q)$ for general $(s,\mf q)$.

\subsection{Continuity of the two-stage SLCP solution}

To compare solutions of the two-stage SLCP at different times,
we need to control the difference between expectations under
$P_t$ and $P_s$. The following lemma provides such an estimate.

\begin{lemma}[Change of measure]\label{lem:tv}
Let $\chi:\Xi\to[0,\infty)$ be measurable.
Suppose there exists a constant $\bar\chi>0$ such that
$
\sup_{r\in[0,T]}\mb E_{P_r}[\chi^2]
\le \bar\chi^2.
$
Then the following statements hold for all $t,s\in[0,T]$.
\begin{enumerate}[(i)]
\item $\displaystyle
\int_\Xi \chi\,d|P_t-P_s|
\le\sqrt2\,\bar \chi\,(L_P|t-s|)^{1/2}$.
\item If $\mf v:\Xi\to\mb R^p$ is measurable and
$\|\mf v(\xi)\|\le a+\chi(\xi)$ for some $a\ge0$, then
$
\|\mb E_{P_t}[\mf v]-\mb E_{P_s}[\mf v]\|
\le aL_P|t-s|+\sqrt2\,\bar \chi\,(L_P|t-s|)^{1/2}.
$
\end{enumerate}
\end{lemma}

\begin{proof}
Since $|P_t-P_s|\le P_t+P_s$,
the Cauchy--Schwarz inequality and
Assumption~\ref{ass:kernel} give
$
\begin{aligned}
\int_\Xi \chi\,d|P_t\!-\!P_s|
\!\le\!\big(\mb E_{P_t}[\chi^2]
\!+\!\mb E_{P_s}[\chi^2]\big)^{1/2}
\|P_t-P_s\|_{\tv}^{1/2}\!\le\!\sqrt{2}\,\bar \chi\,(L_P|t\!-\!s|)^{1/2},
\end{aligned}
$
which proves (i). Part (ii) follows from (i),
$\|\mf v(\xi)\|\le a+\chi(\xi)$,
and Assumption~\ref{ass:kernel}.
\end{proof}

We next derive growth and continuity estimates for the
two-stage SLCP solution. The time variation of $P_t$
is handled using Lemma~\ref{lem:tv}. Define
\begin{equation}\label{eq:omega}
\begin{aligned}
c_y&:=\frac{\bar\kappa_q+\bar\theta}{\underline\kappa},
\qquad
c_u:=(1+c_y)\bar\theta,
\qquad
c_z:=\big(c_y^2+c_u^2\big)^{1/2},\\
\omega(\tau)&:=(c_y+c_u)\tau
+\frac{\sqrt2\,(1+\bar\theta)c_u}{\underline\kappa}
(L_P\tau)^{1/2},\qquad\tau\ge0.
\end{aligned}
\end{equation}

\begin{lemma}[Growth and continuity of the two-stage SLCP solution]\label{lem:est}
For all $t,s,r\in[0,T]$ and $\mf x,\mf x'\in\mb R^n$, the following statements hold.

\noindent(i) $\|\mf y(t,\mf x)\|\le c_y\rho(t,\mf x)$,
$\|\mf z(t,\mf x)\|_{\mc H_{P_r}}\le c_z\rho(t,\mf x)$,
and for every $\xi\in\Xi$,
\begin{equation}\label{eq:u-pointwise}
\|\mf u(t,\mf x,\xi)\|\le(1+c_y)\rho(t,\mf x)\theta(\xi).
\end{equation}

\noindent(ii) $\|\mf z(t,\mf x)-\mf z(t,\mf x')\|_{\mc H_{P_r}}
\le c_z\|\mf x-\mf x'\|$.

\noindent(iii) $
\|\mf z(t,\mf x)-\mf z(s,\mf x)\|_{\mc H_{P_s}}
\le\rho(t,\mf x)\,\omega(|t-s|)$.
\end{lemma}

\begin{proof}
(i) Let $\mf q':=\max\{q(t,\mf x),\mf 0\}$. Since
$\mf 0\in\sol_{P_t}(\mf q')$, Lemma~\ref{lem:wellposed}(iii) with
$\mf q=q(t,\mf x)$, \eqref{eq:q-bounds} and \eqref{eq:g-ratio} give
\[
\|\mf y(t,\mf x)\|\le\frac1{\underline\kappa}\Big(\|f(t,\mf x)\|
+\mb E_{P_t}\Big[\|M(\xi)\|\frac{\|g(t,\mf x,\xi)\|}{\kappa(\xi)}\Big]\Big)
\le c_y\rho(t,\mf x).
\]
For every $\xi\in\Xi$, the first bound in \eqref{eq:u-lip-y}
and Assumption~\ref{ass:q} give
$\kappa(\xi)\|\mf u(t,\mf x,\xi)\|
\le\|M(\xi)\|\|\mf y(t,\mf x)\|+\|g(t,\mf x,\xi)\|
\le(1+c_y)\rho(t,\mf x)(\|M(\xi)\|+\kappa_q(\xi))$,
which is \eqref{eq:u-pointwise}. Assumption~\ref{ass:joint} then gives
$\mb E_{P_r}\|\mf u(t,\mf x,\xi)\|^2\le c_u^2\rho(t,\mf x)^2$, and hence
$\|\mf z(t,\mf x)\|_{\mc H_{P_r}}\le c_z\rho(t,\mf x)$.

\smallskip\noindent
(ii) In the same way, Lemma~\ref{lem:wellposed}(iii) with
$\mf q=q(t,\mf x)$ and $\mf q'=q(t,\mf x')$ gives
$\|\mf y(t,\mf x)-\mf y(t,\mf x')\|\le c_y\|\mf x-\mf x'\|$, and the
second bound in \eqref{eq:u-lip-y} gives
$\|\mf u(t,\mf x,\xi)-\mf u(t,\mf x',\xi)\|
\le(1+c_y)\theta(\xi)\|\mf x-\mf x'\|$, so that
$\mb E_{P_r}\|\mf u(t,\mf x,\xi)-\mf u(t,\mf x',\xi)\|^2
\le c_u^2\|\mf x-\mf x'\|^2$ by Assumption~\ref{ass:joint}, and (ii) follows.

\smallskip\noindent
(iii) Fix $t,s\in[0,T]$ and $\mf x\in\mb R^n$, let $\tau:=|t-s|$,
$\mf z_r=(\mf y_r,\mf u_r):=\mf z(r,\mf x)\in\mc H_{P_s}$ for $r\in\{t,s\}$ (by (i)), and
$\Delta_{t,s}:=\mb E_{P_t}[M_{12}(\xi)\mf u_t(\xi)]
-\mb E_{P_s}[M_{12}(\xi)\mf u_t(\xi)]$.
Since $\mf u_t$ solves \eqref{eq:lcp-xi} for every $\xi\in\Xi$, we have
$\mf z_t\in\sol_{P_s}(q(t,\mf x)+(\Delta_{t,s},\mf 0_k))$, whereas
$\mf z_s\in\sol_{P_s}(q(s,\mf x))$. Hence Lemma~\ref{lem:wellposed}(iii),
\eqref{eq:q-bounds} and \eqref{eq:g-ratio} give
$\|\mf y_t-\mf y_s\|\le c_y\tau+\|\Delta_{t,s}\|/\underline\kappa$,
and the second bound in \eqref{eq:u-lip-y} gives
$\|\mf u_t(\xi)-\mf u_s(\xi)\|\le\theta(\xi)(\|\mf y_t-\mf y_s\|+\tau)$.
By Assumption~\ref{ass:joint} and $(1+\bar\theta)c_y+\bar\theta=c_y+c_u$,
\[
\|\mf z_t-\mf z_s\|_{\mc H_{P_s}}
\le(1+\bar\theta)\|\mf y_t-\mf y_s\|+\bar\theta\tau
\le(c_y+c_u)\tau+\frac{1+\bar\theta}{\underline\kappa}\|\Delta_{t,s}\|.
\]
Finally, by \eqref{eq:u-pointwise}, $\|M_{12}(\xi)\mf u_t(\xi)\|\le \chi(\xi)$
with $\chi(\xi):=(1+c_y)\rho(t,\mf x)\|M(\xi)\|\theta(\xi)$, and
$\sup_{r\in[0,T]}\mb E_{P_r}[\chi^2]\le c_u^2\rho(t,\mf x)^2$ by
Assumption~\ref{ass:joint}. Thus Lemma~\ref{lem:tv}(i) gives
$\|\Delta_{t,s}\|\le\sqrt2\,c_u\rho(t,\mf x)(L_P\tau)^{1/2}$,
and (iii) follows since $\rho(t,\mf x)\ge1$.
\end{proof}

\subsection{Continuity of the vector field}
We now use the growth and continuity estimates for the two-stage SLCP
solution to study the vector field of the differential two-stage SLCP.
Define
\begin{equation}\label{eq:F-def}
F(t,\mf x):=
\mb E_{P_t}[\Phi(t,\mf x,\mf z(t,\mf x,\xi),\xi)],
\qquad
L_F:=\bar\kappa_\Phi(1+c_z).
\end{equation}
The integrand is measurable by Assumption~\ref{ass:Phi}
and Lemma~\ref{lem:wellposed}.

\begin{lemma}[Growth and continuity of $F$]\label{lem:F}
There is a constant $c_F\geq 0$, depending only on
$\underline\kappa,\bar\theta,L_P,\bar\kappa_\Phi,\bar\kappa_q$
and $\bar\kappa_{\Phi,q}$, such that for all $t,s\in[0,T]$
and $\mf x,\mf x'\in\mb R^n$, the following statements hold.
\begin{enumerate}[(i)]
\item $F(t,\mf x)$ is well defined and
$\|F(t,\mf x)\|\le L_F\rho(t,\mf x)$.
\item $\|F(t,\mf x)-F(t,\mf x')\|
\le L_F\|\mf x-\mf x'\|$.
\item $\|F(t,\mf x)-F(s,\mf x)\|
\le c_F\rho(t,\mf x)(|t-s|^{1/2}+|t-s|)$.
\end{enumerate}
\end{lemma}

\begin{proof}
(i) By Assumption~\ref{ass:Phi}, the Cauchy--Schwarz
inequality, and Lemma~\ref{lem:est}(i),
\[
\|F(t,\mf x)\|
\le\mb E_{P_t}\big[
\kappa_\Phi(\xi)
(\rho(t,\mf x)+\|\mf z(t,\mf x,\xi)\|)\big]
\le\bar\kappa_\Phi
\big(\rho(t,\mf x)
+\|\mf z(t,\mf x)\|_{\mc H_{P_t}}\big)\le L_F\rho(t,\mf x).
\]

\smallskip\noindent
(ii) By Assumption~\ref{ass:Phi}, the Cauchy--Schwarz
inequality, and Lemma~\ref{lem:est}(ii),
\[
\|F(t,\mf x)-F(t,\mf x')\|\le\mb E_{P_t}\big[
\kappa_\Phi(\xi)
(\|\mf x-\mf x'\|
+\|\mf z(t,\mf x,\xi)
-\mf z(t,\mf x',\xi)\|)\big]
\le L_F\|\mf x-\mf x'\|.
\]

\noindent(iii) Fix $\mf x$ and write
$\Phi_r(\xi):=\Phi(r,\mf x,\mf z(r,\mf x,\xi),\xi)$ for $r\in\{t,s\}$.
Then
$
F(t,\mf x)-F(s,\mf x)
=\big(\mb E_{P_t}[\Phi_t]-\mb E_{P_s}[\Phi_t]\big)
+\mb E_{P_s}[\Phi_t-\Phi_s].
$
For the first term, Assumption~\ref{ass:Phi} and \eqref{eq:u-pointwise}
give, for every $\xi\in\Xi$,
\[
\|\Phi_t(\xi)\|
\le\kappa_\Phi(\xi)\big(\rho(t,\mf x)+\|\mf z(t,\mf x,\xi)\|\big)
\le\rho(t,\mf x)(1+c_y)\,\kappa_\Phi(\xi)
\big(1+\theta(\xi)\big).
\]
Thus, by Minkowski's inequality and
Assumptions~\ref{ass:Phi} and \ref{ass:joint}, for every $r\in[0,T]$,
$
\mb E_{P_r}\big[\|\Phi_t\|^2\big]^{1/2}
\le\rho(t,\mf x)(1+c_y)
\big(\bar\kappa_\Phi+\bar\kappa_{\Phi,q}\big).
$
Applying Lemma~\ref{lem:tv}(ii) with
$\mf v=\Phi_t$, $a=0$, and $\chi=\|\Phi_t\|$
we obtain
\begin{equation}\label{eq:first}
\|\mb E_{P_t}[\Phi_t]-\mb E_{P_s}[\Phi_t]\|
\le\sqrt2\,\rho(t,\mf x)(1+c_y)
\big(\bar\kappa_\Phi+\bar\kappa_{\Phi,q}\big)
(L_P|t-s|)^{1/2}.
\end{equation}
For the second term, Assumption~\ref{ass:Phi},
Lemma~\ref{lem:est}(iii), and the Cauchy--Schwarz inequality
yield
\begin{align}
\mb E_{P_s}\|\Phi_t-\Phi_s\|
&\le\mb E_{P_s}\big[\kappa_\Phi(\xi)
\big(|t-s|+\|\mf z(t,\mf x,\xi)-\mf z(s,\mf x,\xi)\|\big)\big]\nonumber\\
&\le\bar\kappa_\Phi|t-s|
+\bar\kappa_\Phi
\|\mf z(t,\mf x)-\mf z(s,\mf x)\|_{\mc H_{P_s}}\nonumber\\
&\le\rho(t,\mf x)\bar\kappa_\Phi
\big(|t-s|+\omega(|t-s|)\big).\label{eq:second}
\end{align}
Combining \eqref{eq:omega}, \eqref{eq:first}, and \eqref{eq:second} proves (iii).

\end{proof}

\subsection{Existence and uniqueness of classical solutions}

With the preceding estimates, we are now in a position
to establish the main theorem on the existence and uniqueness
of classical solutions to the differential two-stage SLCP.
\begin{theorem}[Classical solution]\label{thm:classical}
Let $\mf x_0\in\mb R^n$ and $T>0$. Then the differential two-stage SLCP
\eqref{DSVI} has a unique classical solution $(x,z)$ on $[0,T]$
in the sense of Definition~\ref{def:classical}. Uniqueness is
understood pointwise for $x(t)$ and in
$\mc H_{P_t}$ for $z(t)$. Moreover, for all $t\in[0,T]$,
\begin{equation}\label{eq:rho-bound}
\rho(t,x(t))\le(1+\|\mf x_0\|)e^{(1+L_F)t}\le\bar\rho_T,\qquad
\|\dot x(t)\|\le L_F\bar\rho_T,
\end{equation}
where $\bar\rho_T:=(1+\|\mf x_0\|)e^{(1+L_F)T}$.
\end{theorem}

\begin{proof}
By Lemma~\ref{lem:F}, $F$ is continuous in $t$,
globally Lipschitz in $\mf x$ uniformly in $t$, and has
linear growth. Therefore, the initial value problem
$\dot x(t)=F(t,x(t))$, $x(0)=\mf x_0$, has a unique solution
$x\in C^1([0,T];\mb R^n)$ (see, e.g.,
\cite[Corollary~2.6]{teschl2012ordinary}).
Set $z(t):=\mf z(t,x(t))$ and define $y(t)$ and
$u(t,\xi)$ accordingly. Then the two-stage complementarity
conditions and the differential equation hold by construction.

For the bounds, Lemma~\ref{lem:F}(i) gives
$\|\dot x(t)\|\le L_F\rho(t,x(t))$ and
\[
\begin{aligned}
\rho(t,x(t))
\!\le\!1+\|\mf x_0\|+t
+L_F\int_0^t\rho(s,x(s))\,ds\!\le\!1+\|\mf x_0\|
+(1+L_F)\int_0^t\rho(s,x(s))\,ds,
\end{aligned}
\]
where the last inequality uses $\rho(s,x(s))\ge1$.
Gronwall's inequality proves \eqref{eq:rho-bound}.

Next, by Lemma~\ref{lem:est}(ii)--(iii), we have
\[
\begin{aligned}
\|y(t)-y(s)\|
&\le\|\mf z(t,x(t))-\mf z(s,x(t))\|_{\mc H_{P_s}}+\|\mf z(s,x(t))-\mf z(s,x(s))\|_{\mc H_{P_s}}\\
&\le\bar\rho_T
\big(\omega(|t-s|)+c_zL_F|t-s|\big),
\end{aligned}
\]
which yields the continuity of $y$.

For every $\xi\in\Xi$, the second bound in \eqref{eq:u-lip-y}
and Assumption~\ref{ass:q} imply
\[
\begin{aligned}
\|u(t,\xi)-u(s,\xi)\|
\le\frac1{\kappa(\xi)}\big(\|M_{21}(\xi)\|\|y(t)-y(s)\|+\kappa_q(\xi)(|t-s|+\|x(t)-x(s)\|)\big).
\end{aligned}
\]
Thus $u(\cdot,\xi)$ is continuous on $[0,T]$ for every
$\xi\in\Xi$, whereas $u(t,\cdot)$ is measurable by
Lemma~\ref{lem:wellposed}. The chosen representative is therefore
jointly measurable.
This verifies Definition~\ref{def:classical}.

Finally, any other classical solution $(\tilde x,\tilde z)$
satisfies $\tilde z(t)=\mf z(t,\tilde x(t))$ in
$\mc H_{P_t}$ by Lemma~\ref{lem:wellposed}, and hence
$\dot{\tilde x}(t)=F(t,\tilde x(t))$ with
$\tilde x(0)=\mf x_0$. Uniqueness of the ODE solution gives
$\tilde x(t)=x(t)$ and then $\tilde z(t)=z(t)$ in
$\mc H_{P_t}$ for every $t\in[0,T]$.
\end{proof}

\begin{remark}[Joint moment conditions]\label{rem:joint}
Assumption~\ref{ass:joint} cannot be dropped or replaced by separate
second-moment bounds on $\|M\|$, $\kappa_\Phi$ and $\kappa_q$. We give
three examples, with $\Xi=(0,1]$, $\mu$ the Lebesgue measure and
$n=m=k=1$ throughout. Example (a) is static. Examples (b) and (c) use $\zeta(\xi)=\xi^{-1/4}$ and
$P_t=(1-t)P_0+tU_t$ for $0<t\le\frac12$,
where $P_0$ has density $p_0(\xi)=2\xi$ and
$U_t$ is uniform on $(t^2,2t^2]$.
For $t>\frac12$, set $P_t=P_{1/2}$. Then we have
\begin{equation}\label{eq:joint-moments}
\mb E_{P_t}[\zeta^2]
=\frac43(1-t)+2(\sqrt2-1),\quad
\mb E_{P_t}[\zeta^4]=2(1-t)+\frac{\ln2}{t},\quad 0<t\le\frac12,
\quad \mb E_{P_0}[\zeta^2]=\frac43.
\end{equation}
This family satisfies Assumption~\ref{ass:kernel} with $L_P=8$.
Indeed, for any continuous $\varphi:\Xi\to[-1,1]$ and
$0<r<\frac12$,
$|\frac{d}{dr}\int_\Xi\varphi\,p_r\,d\xi|
=|-r^{-2}\int_{r^2}^{2r^2}\varphi\,d\xi
+4\varphi(2r^2)-2\varphi(r^2)
-2\int_\Xi\varphi\xi\,d\xi|\leq8$.
Since $P_t$ is Lipschitz continuous in total variation on $t\in[0,\frac{1}{2}]$ and remains constant for $t\geq\frac12$,
it follows that $\|P_t-P_s\|_{\mathrm{TV}}\leq8|t-s|$
for all $t,s\in[0,T]$.
\begin{enumerate}[(a)]
\item \emph{First condition.} Let $P_t$ be uniform for all $t$,
$M(\xi)=\left(\begin{smallmatrix}1&\xi/2\\-\xi/2&\xi\end{smallmatrix}\right)$,
$f\equiv0$, $g\equiv-1$ and $\Phi=\mf u$.
Assumptions~\ref{ass:Phi}--\ref{ass:kernel} hold with
$\kappa(\xi)=\xi$ and $\kappa_q\equiv1$, but $\mb E_{P_t}[\theta^2]=\infty$.
The unique pointwise solution of \eqref{DSVI-3}--\eqref{DSVI-4} is
$\mf y=\mf 0$, $\mf u(\xi)=1/\xi\notin L^2_{P_t}(\mb R)$, and
$\mb E_{P_t}[\Phi]=\infty$.
\item \emph{Factor $(1+\|M(\xi)\|)^2$.} Let
$M(\xi)=\left(\begin{smallmatrix}1&\zeta(\xi)\\-\zeta(\xi)&1\end{smallmatrix}\right)$,
$f\equiv-c$ with a constant $c>0$, $g(t,\mf x,\xi)=-\zeta(\xi)(1+|\mf x|)$
and $\Phi\equiv0$. Then $\kappa\equiv1$, $\|M\|\le1+\zeta$ and
$\kappa_q=c+\zeta$, so that Assumptions~\ref{ass:Phi}--\ref{ass:kernel},
$\sup_t\mb E_{P_t}[\theta^2]<\infty$ and the second condition in
Assumption~\ref{ass:joint} hold, whereas
$\mb E_{P_t}[(1+\|M\|)^2\theta^2]\ge\mb E_{P_t}[\zeta^4]\to\infty$
as $t\downarrow0$. For $\mf y\ge0$, the second stage gives
$\mf u_{\mf y,\mf g}(\xi)=\zeta(\xi)(\mf y+1+|\mf x|)$, and \eqref{eq:reduced}
yields $\mf y(t,\mf x)=\big(c-a(t)(1+|\mf x|)\big)/(1+a(t))$ with
$a(t):=\mb E_{P_t}[\zeta^2]$ whenever $c>a(t)(1+|\mf x|)$.
Since $\Phi\equiv0$, any classical solution has $x(t)\equiv\mf x_0$, and
by \eqref{eq:joint-moments}, $y(t)=\mf y(t,\mf x_0)$ is discontinuous at
$t=0$ if $c>(\frac43+2(\sqrt2-1))(1+|\mf x_0|)$, contradicting
Definition~\ref{def:classical}.
\item \emph{Second condition.} Let $M(\xi)=I_2$, $f\equiv0$,
$g(t,\mf x,\xi)=-\zeta(\xi)(1+|\mf x|)$ and
$\Phi(t,\mf x,\mf y,\mf u,\xi)=\zeta(\xi)\mf u$.
Then Assumptions~\ref{ass:Phi}--\ref{ass:kernel} and the first condition
in Assumption~\ref{ass:joint} hold with $\kappa\equiv1$ and
$\kappa_\Phi=\kappa_q=\zeta$, while the second condition fails since
$\mb E_{P_t}[\zeta^4]\to\infty$ as $t\downarrow0$.
The unique solution of the two-stage SLCP is
$\mf y(t,\mf x)=\mf 0$ and $\mf u(t,\mf x,\xi)=\zeta(\xi)(1+|\mf x|)$, so that
$F(t,\mf x)=(1+|\mf x|)\mb E_{P_t}[\zeta^2]$ jumps at $t=0$ by
\eqref{eq:joint-moments}. If a classical solution $(x,z)$ existed on
$[0,T]$, continuity of $x$ would give
$\lim_{t\downarrow0}\dot x(t)=(\frac43+2(\sqrt2-1))(1+|\mf x_0|)
\ne\frac43(1+|\mf x_0|)=\dot x(0)$, contradicting continuity of $\dot x$.
Thus no classical solution exists for any $\mf x_0$.
\end{enumerate}
\end{remark}

\section{Discretization and convergence}\label{sec:discretization}

We approximate the differential two-stage SLCP \eqref{DSVI}
by replacing $P_t$ with finite distributions,
solving the resulting two-stage SLCP approximately,
and applying forward Euler time stepping.
The following assumption imposes uniform bounds
and Lipschitz continuity in the scenario variable.
It yields estimates independent of the scenario partition
and the Lipschitz continuity in time required for a first-order Euler bound.

\begin{assumption}[Uniform bounds and Lipschitz continuity in $\xi$]\label{ass:scenario}
The uncertainty set $\Xi\subset\mb R^\ell$ is compact, and the pointwise conditions
in Assumptions~\ref{ass:Phi}, \ref{ass:matrix} and \ref{ass:q} hold for every
$\xi\in\Xi$, with $\kappa_\Phi(\xi)\equiv\bar\kappa_\Phi$,
$\kappa_q(\xi)\equiv\bar\kappa_q$ and $\|M(\xi)\|\le\bar M$.
There exist constants $L_M,L_q,L_\Phi\ge0$
such that, for all $\xi,\xi'\in\Xi$, $t\in[0,T]$,
$\mf x\in\mb R^n$ and $\mf z\in\mb R^{m+k}$,
\begin{align*}
\|M(\xi)-M(\xi')\|&\le L_M\|\xi-\xi'\|,\\
\|q(t,\mf x,\xi)-q(t,\mf x,\xi')\|
&\le L_q\rho(t,\mf x)\|\xi-\xi'\|,\\
\|\Phi(t,\mf x,\mf z,\xi)-\Phi(t,\mf x,\mf z,\xi')\|
&\le L_\Phi\big(\rho(t,\mf x)+\|\mf z\|\big)\|\xi-\xi'\|.
\end{align*}
\end{assumption}

Throughout this section, Assumptions~\ref{ass:kernel} and \ref{ass:scenario}
are in force. We first show that the Lipschitz continuity and positive
definiteness conditions in these assumptions imply
Assumptions~\ref{ass:Phi}--\ref{ass:joint}.
Since $M(\cdot)$ is Lipschitz continuous on the compact set
$\Xi$ and $M(\xi)$ is positive definite for every $\xi\in\Xi$, we have
\begin{equation}\label{eq:eps-uniform}
\varepsilon:=\min\big\{\mf v^\top M(\xi)\mf v:\ \xi\in\Xi,\ \|\mf v\|=1\big\}>0.
\end{equation}
Together with the Lipschitz bounds in Assumptions~\ref{ass:kernel}
and \ref{ass:scenario}, this verifies
Assumptions~\ref{ass:Phi}--\ref{ass:joint} with
$\kappa\equiv\varepsilon$, $\underline\kappa=\varepsilon$,
$\bar\theta=(1+\bar M)(\bar M+\bar\kappa_q)/\varepsilon$, and
$\bar\kappa_{\Phi,q}
=\bar\kappa_\Phi(\bar M+\bar\kappa_q)/\varepsilon$.
The constants $c_y,c_u,c_z,L_F,\bar\rho_T$
in Section~\ref{sec:classical} are specified accordingly.
Consequently, the differential two-stage SLCP admits a unique classical
solution $(x,z)$ on $[0,T]$.

In this uniform setting, we use the Hilbert space structure of
$\mc H_{P_t}$ given by the inner product
$\ip{\mf z}{\mf z'}_{\mc H_{P_t}}:=\mf y^\top\mf y'
+\mb E_{P_t}[\mf u(\xi)^\top\mf u'(\xi)]$, and the closed convex cone
$\mc K_{P_t}:=\mb R^m_+\times L^2_{P_t}(\mb R^k_+)$, which is self-dual
in $\mc H_{P_t}$. Since $\|M(\xi)\|\le\bar M$ on $\Xi$, the block matrix
$M(\xi)$ defines a bounded linear operator $\mc M_t$ on $\mc H_{P_t}$ by
\begin{equation}\label{eq:Mt-uniform}
\begin{aligned}
\mc M_t\mf z
&:=\big(
M_{11}\mf y+\mb E_{P_t}[M_{12}(\xi)\mf u(\xi)],
\ M_{21}(\cdot)\mf y+M_{22}(\cdot)\mf u(\cdot)
\big),\\
\|\mc M_t\mf z\|_{\mc H_{P_t}}&\le\bar M\|\mf z\|_{\mc H_{P_t}},\qquad
\ip{\mc M_t\mf z}{\mf z}_{\mc H_{P_t}}\ge\varepsilon\|\mf z\|_{\mc H_{P_t}}^2,
\end{aligned}
\end{equation}
where the upper bound follows from Jensen's inequality and the lower bound
from \eqref{eq:eps-uniform}. Since $\mc K_{P_t}$ is self-dual, the two-stage SLCP \eqref{eq:tslcp} with $\mf q\in\mc H_{P_t}$ is equivalent to the LCP in $\mc H_{P_t}$: find $\mf z\in\mc H_{P_t}$ such that
\begin{equation}\label{eq:lcp-hilbert}
\mf z\in\mc K_{P_t},\qquad
\mc M_t\mf z+\mf q\in\mc K_{P_t},\qquad
\ip{\mf z}{\mc M_t\mf z+\mf q}_{\mc H_{P_t}}=0,
\end{equation}
which we denote by $\lcp_{\mc H_{P_t}}(\mf q,\mc M_t)$, with solution
set $\sol_{\mc H_{P_t}}(\mf q,\mc M_t)=\sol_{P_t}(\mf q)$. The same applies to the finite-scenario problems introduced below.

\subsection{Scenario discretization by partition}\label{subsec:partition}

Let $\{\Xi_i\}_{i\in[N]}$ be nonempty Borel subsets of $\Xi$ satisfying
$\Xi=\bigcup_{i=1}^N\Xi_i$ and $\Xi_i\cap\Xi_j=\emptyset$ whenever $i\ne j$.
The mesh size of this partition is defined by
\begin{equation}\label{eq:mesh-size}
\delta_N:=\max_{i\in[N]}\operatorname{diam}(\Xi_i),
\qquad
\operatorname{diam}(A):=
\sup_{\xi,\xi'\in A}\|\xi-\xi'\|,\quad A\in\mathcal B(\Xi).
\end{equation}

Choose a representative $\xi_i\in\Xi_i$ for each $i\in[N]$
and set $\pi_i(t):=P_t(\Xi_i)$.
Define
\[
\begin{aligned}
P_t^N(A):=\sum_{i=1}^N\pi_i(t)\mathbf 1_A(\xi_i),\quad
\mu^N(A):=\sum_{i=1}^N\mathbf 1_A(\xi_i),
\end{aligned}
\]
where $\mathbf 1_A$ denotes the indicator function of $A$.
We assume that $\pi_i(t)>0$ for all $t\in[0,T]$ and $i\in[N]$. The representatives are fixed, whereas the probabilities $\pi_i(t)$ may vary
with time. In particular, $\sum_{i=1}^N\pi_i(t)=1$ for every $t\in[0,T]$.
Let $\Xi^N:=\{\xi_1,\dots,\xi_N\}$. Replacing $(\Xi,\mu,P_t)$ by
$(\Xi^N,\mu^N,P_t^N)$ in \eqref{DSVI} gives the $N$-scenario problem
\begin{equation}\label{eq:semi-discrete}
\begin{aligned}
\dot x^N(t)
&=\sum_{i=1}^N\pi_i(t)\Phi\big(t,x^N(t),y^N(t),u_i^N(t),\xi_i\big),
\qquad x^N(0)=\mf x_0,\\
0\le y^N(t)
&\perp M_{11}y^N(t)+\sum_{i=1}^N\pi_i(t)M_{12}(\xi_i)u_i^N(t)
+f(t,x^N(t))\ge0,\\
0\le u_i^N(t)
&\perp M_{21}(\xi_i)y^N(t)+M_{22}(\xi_i)u_i^N(t)
+g(t,x^N(t),\xi_i)\ge0,\quad i\in[N].
\end{aligned}
\end{equation}

We represent the decisions of \eqref{eq:semi-discrete} by arrays
$\mf z=(\mf y,(\mf u_i)_{i=1}^N)$ and write $\mc H_{P_t^N}$ for the space of
arrays with the inner product
$\ip{\mf z}{\mf z'}_{\mc H_{P_t^N}}
:=\mf y^\top\mf y'+\sum_{i=1}^N\pi_i(t)\mf u_i^\top\mf u_i'$,
and $\mc K_{P_t^N}$ for the cone of nonnegative arrays. The operator and right-hand side associated with \eqref{eq:semi-discrete} are
\[
\begin{aligned}
\mc M_t^N\mf z
&:=\Big(M_{11}\mf y+\sum_{i=1}^N\pi_i(t)M_{12}(\xi_i)\mf u_i,\,
\big(M_{21}(\xi_i)\mf y+M_{22}(\xi_i)\mf u_i\big)_{i=1}^N\Big),\\
q^N(t,\mf x)
&:=\Big(f(t,\mf x),\,
\big(g(t,\mf x,\xi_i)\big)_{i=1}^N\Big),
\end{aligned}
\]
and $\mc S^N(t,\mf x)$ denotes the set of arrays satisfying all the
complementarity conditions in \eqref{eq:semi-discrete} at $(t,\mf x)$.
For any array $\mf z=(\mf y,(\mf u_i)_{i=1}^N)$, its piecewise
constant extension is $\hat{\mf z}(\xi):=(\mf y,\mf u_i)$ for
$\xi\in\Xi_i$, so that
$\|\hat{\mf z}\|_{\mc H_{P_t}}=\|\mf z\|_{\mc H_{P_t^N}}$.

The next lemma shows that the finite-scenario system inherits the assumptions
used for the continuous model, with constants independent of the partition.

\begin{lemma}[$N$-scenario well-posedness]\label{lem:inherit}
Under Assumptions~\ref{ass:kernel} and \ref{ass:scenario},
the following statements hold.
\begin{enumerate}[(i)]
\item For all $t,s\in[0,T]$,
$
\|P_t^N-P_s^N\|_{\tv}\le L_P|t-s|.
$
\item For every $(t,\mf x)\in[0,T]\times\mb R^n$,
$\mc S^N(t,\mf x)$ is a singleton.
Moreover, \eqref{eq:semi-discrete} has a unique classical
solution $(x^N,z^N)$ on $[0,T]$, and the estimates of
Section~\ref{sec:classical} hold with the same constants.
In particular,
$
\rho(t,x^N(t))\le\bar\rho_T$ and $
\|\dot x^N(t)\|\le L_F\bar\rho_T.
$
\end{enumerate}
\end{lemma}

\begin{proof}
(i) By Assumption~\ref{ass:kernel},
$
\|P_t^N-P_s^N\|_{\tv}
=\sum_{i=1}^N|\pi_i(t)-\pi_i(s)|\le\sum_{i=1}^N|P_t-P_s|(\Xi_i)
=\|P_t-P_s\|_{\tv}
\le L_P|t-s|.
$

\smallskip\noindent
(ii) By Assumption~\ref{ass:scenario} and \eqref{eq:eps-uniform},
Assumptions~\ref{ass:Phi}--\ref{ass:joint} hold for the $N$-scenario
problem with the same constants, since $P_t^N$ is supported in $\Xi$. In particular, \eqref{eq:Mt-uniform} holds for $\mc M_t^N$ on
$\mc H_{P_t^N}$. The conclusions follow from Lemma~\ref{lem:wellposed}
and Theorem~\ref{thm:classical}.
\end{proof}

We denote the unique solution of the LCPs in \eqref{eq:semi-discrete} at $(t,\mf x)$ by
\begin{equation}\label{eq:solN-notation}
\mf z^N(t,\mf x):=\big(\mf y^N(t,\mf x),(\mf u_i^N(t,\mf x))_{i=1}^N\big),
\end{equation}
with scenario values
$
\mf z^N(t,\mf x,\xi_i):=\big(\mf y^N(t,\mf x),\mf u_i^N(t,\mf x)\big),
\  i\in[N].
$
Denote its piecewise constant extension by
$\hat{\mf z}^N(t,\mf x)$.
The corresponding vector field is
\begin{equation}\label{eq:FN-def}
F^N(t,\mf x):=\sum_{i=1}^N\pi_i(t)
\Phi\big(t,\mf x,\mf z^N(t,\mf x,\xi_i),\xi_i\big).
\end{equation}
Thus the classical solution of \eqref{eq:semi-discrete} satisfies
$\dot x^N(t)=F^N(t,x^N(t))$ and $z^N(t)=\mf z^N(t,x^N(t))$.
We first estimate the error caused by replacing the original data with their
values at the scenario representatives.

\begin{lemma}[Error of the $N$-scenario problem]\label{lem:error}
Let
$C_1:=\varepsilon^{-1}(\sqrt2 L_Mc_z+L_q)$ and
$C_2:=\bar\kappa_\Phi C_1+L_\Phi(1+c_z)$.
Then, for every $(t,\mf x)\in[0,T]\times\mb R^n$, the following statements hold.
\begin{enumerate}[(i)]
    \item $\displaystyle
\|\mf z(t,\mf x)-\hat{\mf z}^N(t,\mf x)\|_{\mc H_{P_t}}
\le C_1\delta_N\rho(t,\mf x)$.
\item $\displaystyle
\|F(t,\mf x)-F^N(t,\mf x)\|\le C_2\delta_N\rho(t,\mf x)$.
\end{enumerate}
\end{lemma}

\begin{proof}
(i) For $\xi\in\Xi_i$, define the piecewise constant data
$\hat M(\xi):=M(\xi_i)$ and
$\hat q(t,\mf x,\xi):=q(t,\mf x,\xi_i)$.
Let $\hat{\mc M}_t$ denote the operator \eqref{eq:Mt-uniform}
with $M$ replaced by $\hat M$.
The solution array satisfies the scenario complementarity
conditions in \eqref{eq:semi-discrete}.
Since $\pi_i(t)=P_t(\Xi_i)$, its piecewise constant
extension satisfies
$\hat{\mf z}^N(t,\mf x)\in
\sol_{\mc H_{P_t}}(\hat q(t,\mf x),\hat{\mc M}_t)$.

Write $\mf z:=\mf z(t,\mf x)$ and
$\hat{\mf z}:=\hat{\mf z}^N(t,\mf x)$.
Combining the variational inequalities associated with
\eqref{eq:lcp-hilbert} for $\mf z$ and $\hat{\mf z}$,
and using \eqref{eq:Mt-uniform}, we obtain
\[
\begin{aligned}
\varepsilon\|\mf z-\hat{\mf z}\|_{\mc H_{P_t}}^2
&\le\ip{(\mc M_t-\hat{\mc M}_t)\hat{\mf z}
+q(t,\mf x)-\hat q(t,\mf x)}
{\hat{\mf z}-\mf z}_{\mc H_{P_t}}\\
&\le\Big(\|(\mc M_t-\hat{\mc M}_t)\hat{\mf z}\|_{\mc H_{P_t}}
+\|q(t,\mf x)-\hat q(t,\mf x)\|_{\mc H_{P_t}}\Big)
\|\mf z-\hat{\mf z}\|_{\mc H_{P_t}}.
\end{aligned}
\]
Consequently,
\begin{equation}\label{eq:consistency-basic}
\begin{aligned}
\|\mf z-\hat{\mf z}\|_{\mc H_{P_t}}
\le\frac1\varepsilon\Big(
&\|(\mc M_t-\hat{\mc M}_t)\hat{\mf z}\|_{\mc H_{P_t}}+\|q(t,\mf x)-\hat q(t,\mf x)\|_{\mc H_{P_t}}\Big).
\end{aligned}
\end{equation}

We now estimate the two terms on the right-hand side.
Let $D(\xi):=M(\xi)-\hat M(\xi)$ and write
$\hat{\mf z}=(\hat{\mf y},\hat{\mf u})$.
By Assumption~\ref{ass:scenario},
$\|D(\xi)\|\le L_M\delta_N$ for every $\xi\in\Xi$,
and $D_{11}=0$ since $M_{11}$ is independent of $\xi$.
Using Jensen's inequality, we obtain
\[
\begin{aligned}
\|(\mc M_t-\hat{\mc M}_t)\hat{\mf z}\|_{\mc H_{P_t}}^2
&=\big\|\mb E_{P_t}[D_{12}(\xi)\hat{\mf u}(\xi)]\big\|^2
+\mb E_{P_t}\|D_{21}(\xi)\hat{\mf y}
+D_{22}(\xi)\hat{\mf u}(\xi)\|^2\\
&\le(L_M\delta_N)^2
\Big(\mb E_{P_t}\|\hat{\mf u}(\xi)\|^2
+\mb E_{P_t}\|\hat{\mf z}(\xi)\|^2\Big)\le2(L_M\delta_N)^2
\|\hat{\mf z}\|_{\mc H_{P_t}}^2.
\end{aligned}
\]
Similarly, Assumption~\ref{ass:scenario} gives
\[
\begin{aligned}
\|q(t,\mf x)-\hat q(t,\mf x)\|_{\mc H_{P_t}}^2
&=\sum_{i=1}^N\int_{\Xi_i}
\|g(t,\mf x,\xi)-g(t,\mf x,\xi_i)\|^2\,dP_t(\xi)\le L_q^2\delta_N^2\rho(t,\mf x)^2.
\end{aligned}
\]
Moreover, by Lemmas~\ref{lem:inherit} and \ref{lem:est}(i),
$
\|\hat{\mf z}\|_{\mc H_{P_t}}
=\|\mf z^N(t,\mf x)\|_{\mc H_{P_t^N}}
\le c_z\rho(t,\mf x).
$
Substituting these estimates into \eqref{eq:consistency-basic},
we obtain
\[
\begin{aligned}
\|\mf z-\hat{\mf z}\|_{\mc H_{P_t}}
&\le\frac1\varepsilon
\big(\sqrt2 L_Mc_z
+L_q\big)\delta_N\rho(t,\mf x)=C_1\delta_N\rho(t,\mf x),
\end{aligned}
\]
which proves (i).

\smallskip\noindent
(ii) Since $\pi_i(t)\!=\!P_t(\Xi_i)$, we can write
$
F^N(t,\mf x)
\!=\!\sum_{i=1}^N\int_{\Xi_i}
\Phi(t,\mf x,\hat{\mf z}^N(t,\mf x,\xi),\xi_i)
\,dP_t(\xi).
$
Comparing with \eqref{eq:F-def}
and adding and subtracting
$\Phi(t,\mf x,\hat{\mf z}^N(t,\mf x,\xi),\xi)$,
we obtain from Assumptions~\ref{ass:Phi}
and \ref{ass:scenario} that
\[
\begin{aligned}
\|F(t,\mf x)-F^N(t,\mf x)\|
&\le\sum_{i=1}^N\int_{\Xi_i}
\bar\kappa_\Phi
\|\mf z(t,\mf x,\xi)-\hat{\mf z}^N(t,\mf x,\xi)\|
\,dP_t(\xi)\\
&\quad+\sum_{i=1}^N\int_{\Xi_i}
L_\Phi\delta_N
\big(\rho(t,\mf x)+\|\hat{\mf z}^N(t,\mf x,\xi)\|\big)
\,dP_t(\xi)\\
&\le\bar\kappa_\Phi
\|\mf z(t,\mf x)-\hat{\mf z}^N(t,\mf x)\|_{\mc H_{P_t}}+L_\Phi\delta_N
\big(\rho(t,\mf x)
+\|\hat{\mf z}^N(t,\mf x)\|_{\mc H_{P_t}}\big),
\end{aligned}
\]
where the last inequality follows from the Cauchy--Schwarz
inequality and \eqref{eq:z}.
Using (i) and Lemma~\ref{lem:est}(i), we obtain
$
\|F(t,\mf x)-F^N(t,\mf x)\|
\le\big(\bar\kappa_\Phi C_1
+L_\Phi(1+c_z)\big)
\delta_N\rho(t,\mf x)=C_2\delta_N\rho(t,\mf x),
$
which proves (ii).
\end{proof}

Lemma~\ref{lem:error} yields the following error bounds
for the classical trajectories.

\begin{theorem}[Convergence of the $N$-scenario classical solutions]
\label{thm:semi-discrete}
Let $(x,z)$ and $(x^N,z^N)$ be the classical solutions of \eqref{DSVI}
and \eqref{eq:semi-discrete} with the same initial state $\mf x_0$.
Set $\hat z^N(t):=\hat{\mf z}^N(t,x^N(t))$ and
$C_3:=C_1+c_zC_2Te^{L_FT}$. Then
\begin{equation}\label{eq:semi-discrete-rate}
\left\{
\begin{aligned}
\sup_{t\in[0,T]}\|x(t)-x^N(t)\|
&\le C_2Te^{L_FT}\bar\rho_T\delta_N,\\
\sup_{t\in[0,T]}\|z(t)-\hat z^N(t)\|_{\mc H_{P_t}}
&\le C_3\bar\rho_T\delta_N.
\end{aligned}
\right.
\end{equation}
Here $\delta_N$ is defined in \eqref{eq:mesh-size},
$C_1,C_2$ in Lemma~\ref{lem:error},
$L_F$ in \eqref{eq:F-def}, and
$\bar\rho_T$ in \eqref{eq:rho-bound}.
\end{theorem}

\begin{proof}
Since $x(0)=x^N(0)=\mf x_0$, we obtain from Lemma~\ref{lem:error}(ii)
and the uniform $L_F$-Lipschitz continuity of $F^N$ that
\[
\begin{aligned}
\|x(t)-x^N(t)\|
&\le\int_0^t\|F(s,x(s))-F^N(s,x(s))\|\,ds+L_F\int_0^t\|x(s)-x^N(s)\|\,ds\\
&\le C_2\bar\rho_T\delta_N t
+L_F\int_0^t\|x(s)-x^N(s)\|\,ds,
\end{aligned}
\]
where the last inequality follows from \eqref{eq:rho-bound}.
Gronwall's inequality then yields
$
\|x(t)-x^N(t)\|
\le C_2\bar\rho_T\delta_N te^{L_Ft}
\le C_2Te^{L_FT}\bar\rho_T\delta_N,
$
proving the first estimate in \eqref{eq:semi-discrete-rate}.

We now estimate the error between
$z(t)=\mf z(t,x(t))$ and
$\hat z^N(t)=\hat{\mf z}^N(t,x^N(t))$.
Inserting the solution $\mf z(t,x^N(t))$
and applying Lemmas~\ref{lem:est}(ii) and \ref{lem:error}(i),
we obtain
\[
\begin{aligned}
&\|z(t)-\hat z^N(t)\|_{\mc H_{P_t}}
\le\|\mf z(t,x(t))-\mf z(t,x^N(t))\|_{\mc H_{P_t}}+\|\mf z(t,x^N(t))
-\hat{\mf z}^N(t,x^N(t))\|_{\mc H_{P_t}}\\
&\le c_z\|x(t)-x^N(t)\|
+C_1\delta_N\rho(t,x^N(t))\le\big(c_zC_2Te^{L_FT}
+C_1\big)\bar\rho_T\delta_N=C_3\bar\rho_T\delta_N,
\end{aligned}
\]
where the last inequality follows from the first bound in
\eqref{eq:semi-discrete-rate} and
$\rho(t,x^N(t))\le\bar\rho_T$ by Lemma~\ref{lem:inherit}.
This proves the second bound in \eqref{eq:semi-discrete-rate}.
\end{proof}

Theorem~\ref{thm:semi-discrete} establishes an $O(\delta_N)$
convergence rate for both the state and decision trajectories.
Since the constants in \eqref{eq:semi-discrete-rate} are independent
of $N$, the partition, and the representatives, the convergence
holds uniformly over these choices as $\delta_N\to0$.

\subsection{Residual stopping criterion and time stepping}\label{subsec:euler}

For fixed $(t,\mf x)$, the $N$-scenario two-stage SLCP
is an affine variational inequality on a polyhedral set.
Progressive hedging converges linearly
\cite[Theorems~2--3]{rockafellar2019solving}.
For its implementation in two-stage SLCPs,
see \cite[Section~3.3]{chen2019discrete}.
To obtain error bounds for all scenario components
independently of the scenario weights,
we use the following blockwise residual-based stopping
criterion (see, e.g., \cite{Pang1987}).
For $(t,\mf x)\in[0,T]\times\mb R^n$,
let
$\tilde{\mf z}
=(\tilde{\mf y},(\tilde{\mf u}_i)_{i=1}^N)
\in\mc H_{P_t^N}$
be an approximate solution of the $N$-scenario two-stage SLCP.
Define its natural residual by
\begin{equation}\label{eq:residual-blocks}
\begin{aligned}
\mf r_{t,\mf x}(\tilde{\mf z}):=\tilde{\mf z}
-\max\big\{\mf 0,\tilde{\mf z}
-\mc M_t^N\tilde{\mf z}
-q^N(t,\mf x)\big\}=:\big(
\mf r^{\mf y}_{t,\mf x}(\tilde{\mf z}),
(\mf r^{\mf u_i}_{t,\mf x}(\tilde{\mf z}))_{i=1}^N
\big).
\end{aligned}
\end{equation}
The maximum is taken componentwise.
We measure the residual in the block maximum norm
$
\|\mf r_{t,\mf x}(\tilde{\mf z})\|_\infty
:=\max\Big\{
\|\mf r^{\mf y}_{t,\mf x}(\tilde{\mf z})\|,
\max_{i\in[N]}
\|\mf r^{\mf u_i}_{t,\mf x}(\tilde{\mf z})\|
\Big\},
$
where each block is equipped with the Euclidean norm.
This norm controls every block independently of the
scenario weights, allowing us to obtain an error bound
in $\mc H_{P_s}$ for every $s\in[0,T]$.

\begin{lemma}[Error bound for approximate LCP solutions]
\label{lem:tolerance}
Let $(t,\mf x)\in[0,T]\times\mb R^n$ and $\eta\ge0$.
Suppose that
$\tilde{\mf z}=(\tilde{\mf y},
(\tilde{\mf u}_i)_{i=1}^N)$ satisfies
$
\|\mf r_{t,\mf x}(\tilde{\mf z})\|_\infty
\le\eta.
$
Set
$c_0:=(1+\bar M)/\varepsilon$ and
$c_{\mathrm{tol}}:=2c_0(1+\bar M/\varepsilon)$.
Then
\[
\max\Big\{
\|\tilde{\mf y}-\mf y^N(t,\mf x)\|,
\max_{i\in[N]}
\|\tilde{\mf u}_i-\mf u_i^N(t,\mf x)\|
\Big\}
\le\frac{c_{\mathrm{tol}}}{\sqrt2}\eta,
\]
and, for every $s\in[0,T]$,
$
\|\hat{\tilde{\mf z}}
-\hat{\mf z}^N(t,\mf x)\|_{\mc H_{P_s}}
\le c_{\mathrm{tol}}\eta.
$
The constant $c_{\mathrm{tol}}$ is independent of
$N$, $\eta$, $t$, and $s$.
\end{lemma}

\begin{proof}
Write $\mf z^N=\mf z^N(t,\mf x)$
and $\mf q=q^N(t,\mf x)$.
Set $\mf w:=\max\{\mf0,\tilde{\mf z}
-\mc M_t^N\tilde{\mf z}-\mf q\}$,
$\mf r:=\tilde{\mf z}-\mf w$, and
$\mf d:=\tilde{\mf z}-\mf z^N$.
In the following, the inner product and norm
are those of $\mc H_{P_t^N}$.
Since $\mf z^N$ solves the $N$-scenario two-stage SLCP,
the projection property of $\mf w$ and the
complementarity conditions for $\mf z^N$ give
$\langle\mf r-\mc M_t^N\tilde{\mf z}-\mf q,
\mf z^N-\mf w\rangle\le0$
and
$\langle\mc M_t^N\mf z^N+\mf q,
\mf z^N-\mf w\rangle\le0$.
Adding these inequalities and using
$\mf z^N-\mf w=\mf r-\mf d$, we obtain
$\langle\mf r-\mc M_t^N\mf d,
\mf r-\mf d\rangle\le0$.
Expanding this inequality and using \eqref{eq:Mt-uniform}
for $\mc M_t^N$ and the Cauchy--Schwarz inequality, we obtain
\[
\varepsilon\|\mf d\|^2
\le\langle\mc M_t^N\mf d,\mf d\rangle
\le\langle\mf r,\mf d\rangle
+\langle\mc M_t^N\mf d,\mf r\rangle
-\|\mf r\|^2
\le(1+\bar M)\|\mf d\|\|\mf r\|.
\]
Moreover, the stopping criterion and
$\sum_i\pi_i(t)=1$ yield
\[
\|\mf r\|_{\mc H_{P_t^N}}^2
=\|\mf r^{\mf y}_{t,\mf x}(\tilde{\mf z})\|^2
+\sum_i\pi_i(t)
\|\mf r^{\mf u_i}_{t,\mf x}(\tilde{\mf z})\|^2
\le2\eta^2.
\]
Consequently,
$\|\tilde{\mf y}-\mf y^N(t,\mf x)\|
\le\|\mf d\|_{\mc H_{P_t^N}}
\le\sqrt2c_0\eta$.

Next, fix $i\in[N]$ and let $\mf u_i'$
be the unique solution of \eqref{eq:lcp-xi}
at $\xi_i$ with $\mf y=\tilde{\mf y}$.
The residual $\mf r^{\mf u_i}_{t,\mf x}(\tilde{\mf z})$
is precisely the natural residual of
this LCP.
Applying the same projection and
complementarity argument in $\mb R^k$,
with matrix $M_{22}(\xi_i)$, gives
$\|\tilde{\mf u}_i-\mf u_i'\|
\le c_0\|\mf r^{\mf u_i}_{t,\mf x}(\tilde{\mf z})\|
\le c_0\eta$.
By the second bound in \eqref{eq:u-lip-y}
and $\kappa(\xi_i)\ge\varepsilon$, we have
\[
\|\mf u_i'-\mf u_i^N(t,\mf x)\|
\le\frac{\bar M}{\varepsilon}
\|\tilde{\mf y}-\mf y^N(t,\mf x)\|
\le\frac{\sqrt2\bar M}{\varepsilon}c_0\eta.
\]
Hence, $\|\tilde{\mf u}_i-\mf u_i^N(t,\mf x)\|
\le c_0(1+\sqrt2\bar M/\varepsilon)\eta
\le c_{\mathrm{tol}}\eta/\sqrt2$.
Together with the bound for $\tilde{\mf y}$,
this proves the first assertion.

Finally, since $\sum_i\pi_i(s)=1$
for every $s\in[0,T]$,
\[
\|\hat{\tilde{\mf z}}
-\hat{\mf z}^N(t,\mf x)\|_{\mc H_{P_s}}^2
=\|\tilde{\mf y}-\mf y^N(t,\mf x)\|^2
+\sum_i\pi_i(s)
\|\tilde{\mf u}_i-\mf u_i^N(t,\mf x)\|^2
\le c_{\mathrm{tol}}^2\eta^2.
\]
The proof is complete.
\end{proof}

We now combine the residual criterion with forward Euler time stepping.
Let $J\in\mb N$, $h:=T/J$, and $t_j:=jh$ for $j=0,\dots,J$.
Fix a tolerance $\eta\ge0$.
The superscript $\nhe$ indicates the dependence of the numerical solution
on the number of scenarios $N$, the time step $h$ and the tolerance $\eta$.
\begin{algorithm}[htbp]
\caption{Scenario--progressive hedging--Euler scheme}
\label{alg:scheme}
\begin{algorithmic}[1]
\REQUIRE Initial state $\mf x_0$, time horizon $T>0$,
number of steps $J\in\mb N$, scenarios $\{\xi_i\}_{i=1}^N$,
weights $\{\pi_i(t)\}_{i=1}^N$, tolerance $\eta\ge0$.
\STATE Set $h:=T/J$, $t_j:=jh$ for $j=0,\dots,J$, and $\mf x_0^{\nhe}:=\mf x_0$.
\FOR{$j=0,\dots,J$}
\STATE Solve the $N$-scenario two-stage SLCP at $(t_j,\mf x_j^{\nhe})$ by
progressive hedging \cite{chen2019discrete,rockafellar2019solving}.
\STATE Obtain
$\tilde{\mf z}_j^{\nhe}=(\tilde{\mf y}_j^{\nhe},(\tilde{\mf u}_{j,i}^{\nhe})_{i=1}^N)$
with $\|\mf r_{t_j,\mf x_j^{\nhe}}(\tilde{\mf z}_j^{\nhe})\|_\infty\le\eta$.
\IF{$j<J$}
\STATE Set
$\tilde F_j^{\nhe}:=\sum_{i=1}^N\pi_i(t_j)
\Phi(t_j,\mf x_j^{\nhe},\tilde{\mf y}_j^{\nhe},\tilde{\mf u}_{j,i}^{\nhe},\xi_i)$.
\STATE Set $\mf x_{j+1}^{\nhe}:=\mf x_j^{\nhe}+h\tilde F_j^{\nhe}$.
\ENDIF
\ENDFOR
\STATE Output $(\mf x_j^{\nhe},\tilde{\mf z}_j^{\nhe})_{j=0}^J$.
\end{algorithmic}
\end{algorithm}

For $\eta>0$, the convergence of progressive hedging
and the continuity of the residual ensure that the prescribed
tolerance is reached in finitely many iterations.
For $\eta=0$, the lower-level problem must be solved exactly,
and finite termination of progressive hedging is not guaranteed.

To compare the output of Algorithm~\ref{alg:scheme}
with the classical solution, we define the interpolants
\begin{equation}\label{eq:interpolants}
x^{\nhe}(t):=\mf x_j^{\nhe}+(t-t_j)\tilde F_j^{\nhe},\qquad
z^{\nhe}(t):=\tilde{\mf z}_j^{\nhe},
\end{equation}
for $t\in[t_j,t_{j+1})$ and $j<J$, with $x^{\nhe}(T):=\mf x_J^{\nhe}$ and $z^{\nhe}(T):=\tilde{\mf z}_J^{\nhe}$. Thus, $x^{\nhe}$ is continuous and piecewise linear,
while $z^{\nhe}$ is piecewise constant. We denote by
$\hat z^{\nhe}(t)$ the piecewise constant extension of the
array $z^{\nhe}(t)$ to $\Xi$.

\subsection{Convergence of the fully discrete scheme}\label{subsec:discrete-error}

We next combine the scenario approximation, time-stepping, and stopping
errors. Under Assumption~\ref{ass:scenario}, the bounded growth moduli and the uniform modulus $\varepsilon$ in \eqref{eq:eps-uniform} allow us to
establish Lipschitz continuity in time for the finite-scenario vector field.

\begin{theorem}[Error estimate]\label{thm:euler}
Let Assumptions~\ref{ass:kernel} and \ref{ass:scenario} hold,
and let $(x,z)$ be the classical solution of \eqref{DSVI}
on $[0,T]$. There exist constants
$C_4,C_5>0$, depending only on $T$, $\|\mf x_0\|$,
$\varepsilon$ in \eqref{eq:eps-uniform} and the constants in
Assumptions~\ref{ass:Phi}--\ref{ass:kernel} and \ref{ass:scenario},
such that, for every partition and set of representatives
satisfying the conditions in Section~\ref{subsec:partition},
every $J\in\mb N$, and every $\eta\ge0$,
the interpolants \eqref{eq:interpolants} satisfy
\[
\left\{
\begin{aligned}
\sup_{t\in[0,T]}\|x(t)-x^{\nhe}(t)\|
&\le C_4(\delta_N+h+\eta),\\
\sup_{t\in[0,T]}\|z(t)-\hat z^{\nhe}(t)\|_{\mc H_{P_t}}
&\le C_5(\delta_N+h+\eta),
\end{aligned}
\right.
\]
where $x^{\nhe}$ and $\hat z^{\nhe}$ are generated by
Algorithm~\ref{alg:scheme}, which first uses a finite-scenario
partition with mesh size $\delta_N$ in \eqref{eq:mesh-size},
then applies the forward Euler method with $h=T/J$, and
finally solves the resulting finite-dimensional problems
by inexact progressive hedging with tolerance $\eta$.
The constants $C_4$ and $C_5$ are independent of
$N$, $h$, $\eta$, the partition, the representatives,
and the weights $\pi_i(\cdot)$.
\end{theorem}

\begin{proof}
Let $(x^N,z^N)$ be the classical solution of \eqref{eq:semi-discrete}
with $x^N(0)=\mf x_0$. For $j=0,\dots,J$, set
$e_j:=\|x^N(t_j)-\mf x_j^{\nhe}\|$ and $\bar e:=\max_{0\le j\le J}e_j$.
By Lemma~\ref{lem:inherit}, the growth and continuity estimates
of Section~\ref{sec:classical} hold for the $N$-scenario problem
with the same constants. In particular, Lemma~\ref{lem:est}(i) and
$\theta(\xi_i)\le(\bar M+\bar\kappa_q)/\varepsilon\le\bar\theta$ give
\begin{equation}\label{eq:zN-bounds}
\begin{aligned}
\|\mf u_i^N(t,\mf x)\|\le c_u\rho(t,\mf x),\quad
\|\mf z^N(t,\mf x,\xi_i)\|\le(c_y+c_u)\rho(t,\mf x),
\quad i\in[N].
\end{aligned}
\end{equation}
Moreover, $\rho(t,x^N(t))\le\bar\rho_T$ and
$\|\dot x^N(t)\|\le L_F\bar\rho_T$ for all $t\in[0,T]$.

\smallskip\noindent
\emph{Step 1: Lipschitz continuity of $\mf z^N$ in time.}
Fix $t,s\in[0,T]$ and $\mf x\in\mb R^n$, and let
$
\Delta_{t,s}:=\sum_{i=1}^N
(\pi_i(t)-\pi_i(s))M_{12}(\xi_i)\mf u_i^N(t,\mf x)
$
be the measure-change term of the proof of Lemma~\ref{lem:est}(iii)
for the $N$-scenario problem.
We compare $\mf z^N(t,\mf x)$ and $\mf z^N(s,\mf x)$
in $\mc H_{P_s^N}$.
As in the proof of Lemma~\ref{lem:est}(iii),
$\mf z^N(t,\mf x)$ solves the LCP in $\mc H_{P_s^N}$
with operator $\mc M_s^N$ and right-hand side
$q^N(t,\mf x)+(\Delta_{t,s},(\mf 0_k)_{i=1}^N)$.
Comparing its variational inequality with that for
$\mf z^N(s,\mf x)$ and using \eqref{eq:Mt-uniform}, we obtain
$
\|\mf z^N(t,\mf x)-\mf z^N(s,\mf x)\|_{\mc H_{P_s^N}}\le\frac1\varepsilon
\big(\|\Delta_{t,s}\|+\bar\kappa_q|t-s|\big).
$
By Lemma~\ref{lem:inherit}(i) and \eqref{eq:zN-bounds},
$
\|\Delta_{t,s}\|\le\bar M c_u\rho(t,\mf x)
\sum_{i=1}^N|\pi_i(t)-\pi_i(s)|\le\bar M c_uL_P\rho(t,\mf x)|t-s|.
$
Using $\rho(t,\mf x)\ge1$, we obtain
\begin{equation}\label{eq:zN-time-lip}
\|\mf z^N(t,\mf x)-\mf z^N(s,\mf x)\|_{\mc H_{P_s^N}}
\le\rho(t,\mf x)\,\omega_1|t-s|,
\qquad
\omega_1:=\frac{\bar\kappa_q+\bar Mc_uL_P}{\varepsilon}.
\end{equation}

\noindent
\emph{Step 2: Lipschitz continuity of $F^N$ and $\dot x^N$ in time.}
By \eqref{eq:FN-def},
\[
\resizebox{0.94\textwidth}{!}{$\displaystyle
\begin{aligned}
F^N(t,\mf x)-F^N(s,\mf x)
={}&\sum_{i=1}^N(\pi_i(t)-\pi_i(s))
\Phi(t,\mf x,\mf z^N(t,\mf x,\xi_i),\xi_i)\\
&+\sum_{i=1}^N\pi_i(s)
\big(\Phi(t,\mf x,\mf z^N(t,\mf x,\xi_i),\xi_i)
-\Phi(s,\mf x,\mf z^N(s,\mf x,\xi_i),\xi_i)\big).
\end{aligned}
$}
\]
For the first sum, Assumption~\ref{ass:Phi},
\eqref{eq:zN-bounds}, and Lemma~\ref{lem:inherit}(i)
give a bound proportional to $L_P|t-s|$.
For the second sum, we use Assumption~\ref{ass:Phi},
the Cauchy--Schwarz inequality, and \eqref{eq:zN-time-lip}.
Combining these estimates yields
\[
\begin{aligned}
\|F^N(t,\mf x)&-F^N(s,\mf x)\|
\le\bar\kappa_\Phi
\big(1+c_y+c_u\big)
\rho(t,\mf x)L_P|t-s|\\
&\quad+\bar\kappa_\Phi
\Big(|t-s|+
\|\mf z^N(t,\mf x)-\mf z^N(s,\mf x)\|_{\mc H_{P_s^N}}\Big)\le C_6\rho(t,\mf x)|t-s|,
\end{aligned}
\]
where
$
C_6:=\bar\kappa_\Phi
\big(1+(1+c_y+c_u)L_P+\omega_1\big).
$
Since $F^N$ is uniformly $L_F$-Lipschitz in $\mf x$
by Lemmas~\ref{lem:inherit} and \ref{lem:F}(ii), and
$\|x^N(t)-x^N(s)\|\le L_F\bar\rho_T|t-s|$, we have
\[
\begin{aligned}
\|\dot x^N(t)\!-\!\dot x^N(s)\|
&\le\|F^N(t,x^N(t))\!-\!F^N(s,x^N(t))\|+\|F^N(s,x^N(t))\!-\!F^N(s,x^N(s))\|\\
&\le(C_6+L_F^2)\bar\rho_T|t-s|
=:\Lambda_T|t-s|.
\end{aligned}
\]
\noindent
\emph{Step 3: Error of the state approximation.}
At each grid point $t_j$, Lemma~\ref{lem:tolerance}
and Assumption~\ref{ass:Phi} imply
\begin{equation}\label{eq:Ftilde-error}
\begin{aligned}
\|F^N(t_j,\mf x_j^{\nhe})-\tilde F_j^{\nhe}\|
\le\bar\kappa_\Phi
\|\mf z^N(t_j,\mf x_j^{\nhe})-\tilde{\mf z}_j^{\nhe}\|_{\mc H_{P_{t_j}^N}}\le\bar\kappa_\Phi c_{\mathrm{tol}}\eta.
\end{aligned}
\end{equation}
By the Euler update in Algorithm~\ref{alg:scheme}
and the integral equation for $x^N$, we have
\[
\begin{aligned}
x^N(t_{j+1})-\mf x_{j+1}^{\nhe}
={}&x^N(t_j)-\mf x_j^{\nhe}+h\big(F^N(t_j,x^N(t_j))-F^N(t_j,\mf x_j^{\nhe})\big)\\
&+h\big(F^N(t_j,\mf x_j^{\nhe})-\tilde F_j^{\nhe}\big)+\int_{t_j}^{t_{j+1}}
\big(\dot x^N(r)-\dot x^N(t_j)\big)\,dr.
\end{aligned}
\]
The last integral is bounded by $\Lambda_T h^2/2$
from Step~2. Using the state Lipschitz bound for $F^N$
and \eqref{eq:Ftilde-error}, we obtain
\[
e_{j+1}\le(1+hL_F)e_j+\frac12\Lambda_T h^2
+h\bar\kappa_\Phi c_{\mathrm{tol}}\eta,
\qquad j=0,\dots,J-1.
\]
Since $e_0=0$ by $x^N(0)=\mf x_0^{\nhe}=\mf x_0$,
and $(1+hL_F)^j\le e^{L_FT}$,
summing this recursion gives
\begin{equation}\label{eq:E-bound}
\bar e\le Te^{L_FT}
\big(\tfrac12\Lambda_T h+\bar\kappa_\Phi c_{\mathrm{tol}}\eta\big).
\end{equation}
For $t\in[t_j,t_{j+1})$, applying the same estimate
to the interpolant \eqref{eq:interpolants} yields
\[
\|x^N(t)-x^{\nhe}(t)\|
\le(1+hL_F)\bar e+\tfrac12\Lambda_T h^2
+h\bar\kappa_\Phi c_{\mathrm{tol}}\eta.
\]
The inequality also holds at $t=T$, since
$\|x^N(T)-x^{\nhe}(T)\|=e_J\le\bar e$.
Combining it with the first estimate in
\eqref{eq:semi-discrete-rate}, we obtain
\begin{equation}\label{eq:x-total}
\begin{aligned}
\sup_{t\in[0,T]}\|x(t)-x^{\nhe}(t)\|
&\le C_2Te^{L_FT}\bar\rho_T\delta_N
+(1+hL_F)\bar e+\tfrac12\Lambda_T h^2
+h\bar\kappa_\Phi c_{\mathrm{tol}}\eta.
\end{aligned}
\end{equation}

\noindent
\emph{Step 4: Error of the decision approximation.}
By Lemma~\ref{lem:est}(ii), applied to the $N$-scenario
problem by Lemma~\ref{lem:inherit}(ii), and
$\|\hat{\mf z}\|_{\mc H_{P_s}}=\|\mf z\|_{\mc H_{P_s^N}}$,
we have, for all $r,s\in[0,T]$ and $\mf x,\mf x'\in\mb R^n$,
\begin{equation}\label{eq:zN-cross-state}
\|\hat{\mf z}^N(r,\mf x)-\hat{\mf z}^N(r,\mf x')\|_{\mc H_{P_s}}
\le c_z\|\mf x-\mf x'\|.
\end{equation}

Now let $t\in[t_j,t_{j+1})$, or $t=T$ with $j=J$.
Recall that $\hat z^N(t)=\hat{\mf z}^N(t,x^N(t))$.
By inserting the decision $\hat z^N(t_j)$
and the solution $\hat{\mf z}^N(t_j,\mf x_j^{\nhe})$,
the triangle inequality gives
\[
\begin{aligned}
\|z(t)-\hat z^{\nhe}(t)\|_{\mc H_{P_t}}
&\le\|z(t)-\hat z^N(t)\|_{\mc H_{P_t}}+\|\hat z^N(t)-\hat z^N(t_j)\|_{\mc H_{P_t}}\\
&\quad+\|\hat{\mf z}^N(t_j,x^N(t_j))
-\hat{\mf z}^N(t_j,\mf x_j^{\nhe})\|_{\mc H_{P_t}}\\
&\quad+\|\hat{\mf z}^N(t_j,\mf x_j^{\nhe})
-\hat z^{\nhe}(t)\|_{\mc H_{P_t}}.
\end{aligned}
\]
The first term is bounded by $C_3\bar\rho_T\delta_N$
by \eqref{eq:semi-discrete-rate}.

For the second term, inserting
$\hat{\mf z}^N(t,x^N(t_j))$ and applying
\eqref{eq:zN-cross-state} at time $t$,
together with \eqref{eq:zN-time-lip}
at $(t_j,t,x^N(t_j))$, gives
\[
\|\hat z^N(t)-\hat z^N(t_j)\|_{\mc H_{P_t}}
\le c_z\|x^N(t)-x^N(t_j)\|
+\rho(t_j,x^N(t_j))\omega_1|t-t_j|
\le\bar\rho_T(c_zL_F+\omega_1)h.
\]
The third term is at most
$c_ze_j\le c_z\bar e$
by \eqref{eq:zN-cross-state}.
The fourth term is at most $c_{\mathrm{tol}}\eta$
by Lemma~\ref{lem:tolerance} with $s=t$.
Combining these bounds, we obtain
\begin{equation}\label{eq:z-total}
\begin{aligned}
\sup_{t\in[0,T]}\|z(t)-\hat z^{\nhe}(t)\|_{\mc H_{P_t}}
&\le C_3\bar\rho_T\delta_N+\bar\rho_T(c_zL_F+\omega_1)h+c_z\bar e+c_{\mathrm{tol}}\eta.
\end{aligned}
\end{equation}

Substituting \eqref{eq:E-bound} into \eqref{eq:x-total}
and \eqref{eq:z-total}, and using $h\le T$, proves
both estimates with the stated dependence of
$C_4$ and $C_5$.
\end{proof}

\section{Numerical experiments}\label{sec:minimax-experiments}

We consider a two-stage stochastic zero-sum game
between a decision maker and a strategic opponent.
Player 1 minimizes an expected cost, while Player 2
maximizes the same quantity.
Both players choose first-stage actions before
the uncertainty is observed and make recourse decisions
after its realization.
We use a convex-concave quadratic instance of the
two-stage stochastic minimax model in \cite{SunChen2025}
and apply Algorithm~\ref{alg:scheme} to its strategy adjustment dynamics.

\smallskip\noindent{\bf Two-stage game formulation.}
Let $\mf x_1$ and $\mf y_1$ denote the first-stage
actions of Players 1 and 2, respectively.
For each fixed $t$, the game under $P_t$ is
\begin{equation}\label{eq:sc-minimax}
\left\{
\begin{aligned}
&\min_{\mf x_1\in X_1}\ \max_{\mf y_1\in Y_1}
\big\{F_1(\mf x_1,\mf y_1)
+\mb E_{P_t}[\psi_2(\mf x_1,\mf y_1,\xi)]\big\},\\
&\psi_2(\mf x_1,\mf y_1,\xi)
:=\min_{\mf x_2\in X_2(\mf x_1,\xi)}
\max_{\mf y_2\in Y_2(\mf y_1,\xi)}
F_2(\mf x_2,\mf y_2,\xi).
\end{aligned}
\right.
\end{equation}
The recourse sets are
$X_2(\mf x_1,\xi)
=\{\mf x_2\in\mb R^{n_2}:
G(\xi)\mf x_1+W\mf x_2\le b(\xi)\}$
and
$Y_2(\mf y_1,\xi)
=\{\mf y_2\in\mb R^{m_2}:
A(\xi)\mf y_1+B\mf y_2\le c(\xi)\}$.
Here $k_x$ and $k_y$ denote the numbers of
recourse constraints for Players 1 and 2, respectively.
We take $X_1=\mb R^{n_1}$ and
$Y_1=[\mf 0_{m_1},\mf c_1]$ with
$\mf c_1>\mf 0_{m_1}$.
The stage costs of Player 1, equivalently the stage
payoffs of Player 2, are
\begin{equation*}
\left\{
\begin{aligned}
F_1(\mf x_1,\mf y_1)
&\!=\!\tfrac12\mf x_1^\top Q_1\mf x_1+d_1^\top\mf x_1
+\mf x_1^\top O_1\mf y_1-\tfrac12\mf y_1^\top S_1\mf y_1-\beta_1^\top\mf y_1,\\
F_2(\mf x_2,\mf y_2,\xi)
&\!=\!\tfrac12\mf x_2^\top Q_2(\xi)\mf x_2\!+\!d_2(\xi)^\top\mf x_2\!+\!\mf x_2^\top O_2(\xi)\mf y_2
\!-\!\tfrac12\mf y_2^\top S_2(\xi)\mf y_2\!-\!\beta_2(\xi)^\top\mf y_2.
\end{aligned}
\right.
\end{equation*}
The matrices $Q_1$, $S_1$, $Q_2(\xi)$, and $S_2(\xi)$
are symmetric positive definite.
Thus, each stage payoff is strongly convex in
the minimizing player's action and strongly concave
in the maximizing player's action.
For the feasible recourse structures used below,
the game admits a saddle point characterized
by its KKT conditions.

\smallskip\noindent{\bf Differential two-stage SLCP formulation.}
We construct the differential two-stage SLCP using
the KKT conditions of \eqref{eq:sc-minimax}
at each fixed $t$ \cite{SunChen2025}.
Let $\nu_x(\xi)$ and $\nu_y(\xi)$ denote the recourse multipliers,
and let $\lambda_1$ be the multiplier associated with
$\mf y_1\le\mf c_1$.
Set $\mf y=(\mf y_1,\lambda_1)$ and
$\mf u=(\nu_x,\nu_y)$, and define
\begin{equation}\label{eq:elimination-matrices}
\Gamma(\xi):=\begin{pmatrix}Q_2(\xi)&O_2(\xi)\\-O_2(\xi)^\top&S_2(\xi)\end{pmatrix},\ 
R:=\begin{pmatrix}W&0\\0&B\end{pmatrix},\ 
b_2(\xi):=
(d_2(\xi)^\top,\beta_2(\xi)^\top)^\top.
\end{equation}
The second-stage stationarity conditions are
$
Q_2(\xi)\mf x_2+O_2(\xi)\mf y_2+d_2(\xi)
  +W^\top\nu_x(\xi)=0$ and $
-O_2(\xi)^\top\mf x_2+S_2(\xi)\mf y_2+\beta_2(\xi)
  +B^\top\nu_y(\xi)=0.
$
Since the symmetric part of $\Gamma(\xi)$ is positive definite,
$\Gamma(\xi)$ is invertible, and these conditions give
\begin{equation}\label{eq:primal-recovery}
\begin{pmatrix}
\mf x_2(\xi)\\
\mf y_2(\xi)
\end{pmatrix}
=-\Gamma(\xi)^{-1}
\big(b_2(\xi)+R^\top\mf u(\xi)\big).
\end{equation}
Substituting \eqref{eq:primal-recovery} into the KKT conditions
yields \eqref{DSVI-3}--\eqref{DSVI-4} with
\begin{equation*}
\begin{gathered}
M_{11}:=\begin{pmatrix}S_1&I_{m_1}\\-I_{m_1}&0\end{pmatrix},\quad
M_{12}(\xi):=\begin{pmatrix}0&A(\xi)^\top\\0&0\end{pmatrix},\quad
M_{21}(\xi):=-M_{12}(\xi)^\top,\\
M_{22}(\xi):=R\Gamma(\xi)^{-1}R^\top,\qquad
f(\mf x):=(\beta_1-O_1^\top\mf x,\mf c_1),\\
g(\mf x,\xi):=R\Gamma(\xi)^{-1}b_2(\xi)
+(b(\xi)-G(\xi)\mf x,c(\xi)).
\end{gathered}
\end{equation*}
Since $X_1=\mb R^{n_1}$, the first-stage optimality
condition of Player 1 is
$Q_1\mf x_1+O_1\mf y_1+d_1
+\mb E_{P_t}[G(\xi)^\top\nu_x(\xi)]=0$.
We use its negative residual to define the state
vector field
\begin{equation}\label{eq:minimax-phi}
\Phi(\mf x,\mf y,\mf u,\xi)
:=-Q_1\mf x-O_1\mf y_1-d_1-G(\xi)^\top\nu_x.
\end{equation}
The state $x(t)$ describes the adjustment of Player 1's
first-stage action.
At each time, the two-stage complementarity system
determines Player 2's first-stage action and the
constraint multipliers.
Here $t$ is the adjustment time, distinct from
the two stages of the game.

\smallskip\noindent{\bf Differential two-stage SLCP equilibria.}
When $P_t\equiv P$, the system is autonomous, and we
omit the time argument in $\mc S(\mf x)$.
Following \cite{ChenShen2022}, a pair
$(\mf x^e,\mf z^e)\in\mb R^n\times\mc H_P$
is an equilibrium of \eqref{DSVI} if
\begin{equation}\label{eq:equilibrium}
\mb E_P[\Phi(\mf x^e,\mf z^e(\xi),\xi)]=\mf 0_n,
\qquad
\mf z^e\in\mc S(\mf x^e).
\end{equation}
For the unregularized system, \eqref{eq:equilibrium}
and \eqref{eq:primal-recovery} give the KKT conditions
of the game and hence characterize its saddle points.
The operator $\mc M_t$ is monotone, since
for any $\mf v=(\mf y,\mf u)\in\mb R^{m+k}$,
\begin{equation}\label{eq:minimax-monotonicity}
\mf v^\top M(\xi)\mf v
=\mf y_1^\top S_1\mf y_1
+\mf u^\top R\Gamma(\xi)^{-1}R^\top\mf u
\ge0.
\end{equation}
It is not strongly monotone, because the quadratic form
vanishes in the $\lambda_1$ directions.

\smallskip\noindent{\bf Numerical setting.}
All computations were performed in MATLAB R2025a on a 64-bit
Windows system with an Intel Core i9-14900HX processor
(2.20 GHz) and 16 GB RAM.
Following \cite[Section~5.5]{SunChen2025},
we generate numerical game instances with dimensions
$(n_1,m_1,n_2,m_2,k_x,k_y)=(3,2,4,3,2,2)$,
so that $n=3$ and $m=k=4$.
The first-stage parameters are $Q_1=I_3$, $S_1=I_2$,
and $\mf c_1=2\mf 1_2$, while the second-stage constraint
matrices are initialized as
$W_0=[\,I_2\ \ 0_{2\times2}\,]$ and
$B=[\,I_2\ \ 0_{2\times1}\,]$.
The entries of $O_1$ are drawn uniformly from
$[-0.2,0.2]$, and those of $d_1$ and $\beta_1$
from $[-0.5,0.5]$. The data remain fixed across refinements,
except for the perturbations specified below.

We take $\Xi=[-1,1]$.
Each coefficient in
$\{Z_x,Z_y,O_2,G,A,d_2,\beta_2,b,c\}$
is the linear interpolation of two independently generated
endpoint arrays; for example,
$G(\xi)=\frac{1}{2}(G_++G_-)+\frac{\xi}{2}(G_+-G_-)$.
The entries of the endpoint arrays are sampled uniformly
from $[-1,1]$ for $Z_x,Z_y$, $[-0.1,0.1]$ for $O_2$,
$[-0.2,0.2]$ for $G,A$, $[-0.5,0.5]$ for
$d_2,\beta_2$, and $[0.1,0.3]$ for $b,c$.
Here $Z_x(\xi)\in\mb R^{4\times4}$ and
$Z_y(\xi)\in\mb R^{3\times3}$ are auxiliary matrices
used to construct
$
Q_2(\xi)=I_4+0.1Z_x(\xi)^\top Z_x(\xi)$ and $
S_2(\xi)=I_3+0.1Z_y(\xi)^\top Z_y(\xi).
$
Consequently, $Q_2(\xi)\succeq I_4$ and
$S_2(\xi)\succeq I_3$ for every $\xi\in[-1,1]$.
All coefficient functions are bounded and Lipschitz
continuous on $[-1,1]$.

All refinement and regularization tests use $T=2$ and
$P_t(d\xi)=p_t(\xi)\,d\xi$ on $[0,T]$, where
\begin{equation}\label{eq:num-time-dependent-density}
p_t(\xi)=\tfrac12[1+a(t)\xi],\qquad
a(t)=t/T-\tfrac12,\qquad \xi\in[-1,1].
\end{equation}
Since $1/4\le p_t\le3/4$, $\int_{-1}^1p_t\,d\xi=1$, and
$\|P_t-P_s\|_{\tv}
=\int_{-1}^1|p_t-p_s|\,d\xi=|t-s|/(2T)$,
Assumption~\ref{ass:kernel} holds with $L_P=1/4$.
For the partition approximation, we use $N$ equal subintervals
$\Xi_i$, left endpoints $\xi_i=-1+2(i-1)/N$, and $\delta_N=2/N$.
The discrete distribution is
\begin{equation}\label{eq:num-time-dependent-weights}
\left\{
\begin{aligned}
\pi_i(t)&=P_t(\Xi_i)=\frac1N[1+a(t)(\xi_i+1/N)],\\
P_t^N(A)&=\sum_{i=1}^N\pi_i(t)\mathbf 1_A(\xi_i),
\qquad A\in\mathcal B(\Xi).
\end{aligned}
\right.
\end{equation}
The weights are updated at each Euler node.
We apply Algorithm~\ref{alg:scheme} using progressive
hedging (PH) \cite{rockafellar2019solving} with penalty
parameter $r=1$ and at most 3000 iterations per time node.
The eight-dimensional scenario subproblems are solved by an
active-set method with residual tolerance $2\times10^{-12}$.
PH starts from zero and then uses the preceding time node,
except in the tolerance experiments, where it restarts from zero.
We use the natural residual \eqref{eq:residual-blocks}
in the block maximum norm with stopping tolerance $\eta$
for the matrix $M^\alpha(\xi)=M(\xi)+\alpha I_{m+k}$, where $\alpha\ge0$.

We measure the state and decision errors
relative to finer numerical reference solutions:
\begin{equation}\label{eq:num-errors}
\begin{aligned}
E_x:=\sup_{t\in[0,T]}
\|x^{\nhe}(t)-x^{\mathrm{ref}}(t)\|,\quad
E_z:=\sup_{t\in[0,T]}
\|\hat z^{\nhe}(t)-\hat z^{\mathrm{ref}}(t)\|_{\mc H_{P_t}}.
\end{aligned}
\end{equation}
For different scenario sizes, the decisions are represented
on a common refinement of the scenario partitions.
Both left and right decision limits at the time-grid points
are evaluated under the current $P_t$ when computing $E_z$.

\smallskip\noindent{\bf Example 1: duplicated constraints and regularization.}
We replace the second row of $W_0$, $G(\xi)$, and
$b(\xi)$ with their first row, producing duplicated
recourse constraints.
Let $\sigma(\xi):=\nu_{x,1}(\xi)+\nu_{x,2}(\xi)$.
If $\sigma(\xi)>0$ on a set of positive $P_t$-measure,
the splittings
$(\nu_{x,1},\nu_{x,2})
=(\vartheta\sigma,(1-\vartheta)\sigma)$,
$\vartheta\in[0,1]$,
give distinct multipliers with the same primal actions
and stationarity conditions.
To obtain a strongly monotone problem, we introduce the
regularization
\begin{equation}\label{eq:num-regularization}
M^\alpha(\xi):=M(\xi)+\alpha I_{m+k},\qquad
\mc M_t^\alpha:=\mc M_t+\alpha I_{\mc H_{P_t}},\quad\alpha>0.
\end{equation}
By \eqref{eq:minimax-monotonicity}, the regularized operator
is strongly monotone with modulus at least $\alpha$.
Thus, Theorems~\ref{thm:semi-discrete} and \ref{thm:euler}
apply to the regularized differential two-stage SLCP for each fixed $\alpha>0$,
with constants that may depend on $\alpha$.

\begin{figure}[htbp]
\centering
\includegraphics[width=0.9\linewidth]{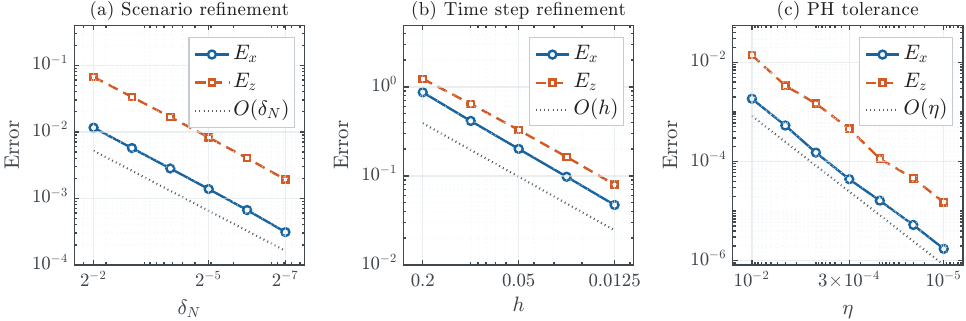}
\caption{Example 1 ($\alpha=10^{-2}$):
State and decision errors as
(a) $\delta_N$ decreases, (b) $h$ decreases, and
(c) $\eta$ decreases. Dotted lines indicate first-order convergence.}
\label{fig:example1-errors}
\end{figure}

We first examine the convergence rates with
$\alpha=10^{-2}$ and $\mf x_0=(12,-12,-12)^\top$.
In Figure~\ref{fig:example1-errors}(a), we vary
$N\in\{8,16,32,64,128,256\}$, with
$(h,\eta)=(0.025,10^{-10})$ and reference $N=2048$.
In Figure~\ref{fig:example1-errors}(b), we successively halve $h$
from $0.2$ to $0.0125$, with $(N,\eta)=(64,10^{-10})$
and reference $h=0.00078125$.
In Figure~\ref{fig:example1-errors}(c), we use
$\eta\in\{10^{-2},3\times10^{-3},10^{-3},3\times10^{-4},
10^{-4},3\times10^{-5},10^{-5}\}$, with $(N,h)=(64,0.025)$
and reference $\eta=10^{-10}$.
The errors decrease at approximately first-order rates
as $\delta_N$, $h$, and $\eta$ decrease, in agreement with
Theorems~\ref{thm:semi-discrete} and \ref{thm:euler}.

We repeat the refinement tests on 50 perturbed instances.
For each instance, we add independent uniform perturbations
from $[-0.03,0.03]$ to $d_1,\beta_1$ and from
$[-0.02,0.02]$ to both endpoint arrays of
$d_2,\beta_2,b,c$. We then replace the second row of
$b$ with its first row to preserve the duplicated constraints.
The same perturbed data are used at all refinement levels.
We retain $\alpha$ and $\mf x_0$ for Table~\ref{tab:observed-orders}.
For scenario refinement, we use
$N\in\{8,16,32,64\}$, $N_{\rm ref}=512$, and
$(h,\eta)=(0.025,10^{-10})$.
For time refinement, we use
$h\in\{0.2,0.1,0.05,0.025\}$,
$h_{\rm ref}=0.0015625$, and $(N,\eta)=(32,10^{-10})$.
For tolerance refinement, we use
$\eta\in\{10^{-2},3\times10^{-3},10^{-3},3\times10^{-4}\}$,
$\eta_{\rm ref}=10^{-10}$, and $(N,h)=(32,0.05)$.
The observed orders $p_x$ and $p_z$ are obtained by
least-squares fitting of $\log E_x$ and $\log E_z$
against $\log s$, respectively, where
$s\in\{\delta_N,h,\eta\}$.
Table~\ref{tab:observed-orders} reports their means
and sample standard deviations over the 50 instances. These results further support the convergence
rates established in Theorems~\ref{thm:semi-discrete}
and \ref{thm:euler}.

\begin{table}[htbp]
\centering\small
\caption{Observed convergence orders in Example 1 over 50 runs
(mean $\pm$ sample standard deviation).}
\label{tab:observed-orders}
\begin{tabular}{@{}lccc@{}}
\toprule
 & $\delta_N$ & $h$ & $\eta$\\
\midrule
$p_x$ & $1.0636\pm2.54\times10^{-4}$
      & $1.0638\pm1.76\times10^{-4}$
      & $1.0700\pm9.93\times10^{-3}$\\
$p_z$ & $1.0381\pm1.89\times10^{-4}$
      & $0.9857\pm1.52\times10^{-6}$
      & $0.9572\pm2.00\times10^{-2}$\\
\bottomrule
\end{tabular}
\end{table}

\begin{figure}[htbp]
\centering
\includegraphics[width=0.9\linewidth]{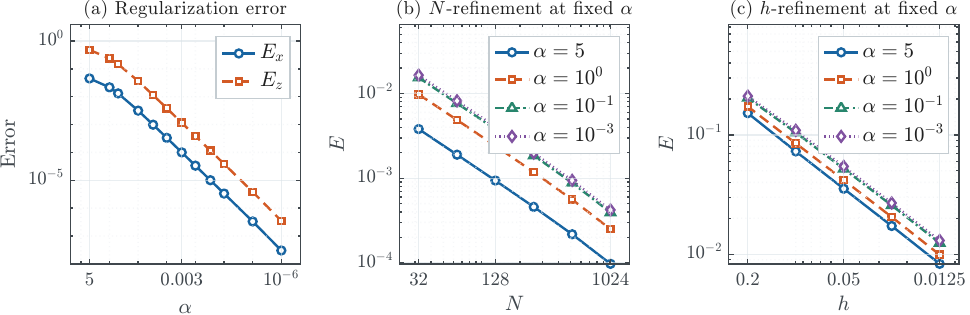}
\caption{Example 1: (a) State and decision errors as $\alpha$
decreases, $\alpha_{\rm ref}=10^{-7}$;
(b) combined error $E$ as $N$ increases for different $\alpha$;
(c) combined error $E$ as $h$ decreases for different $\alpha$.}
\label{fig:regularization}
\end{figure}

We next examine the effects of regularization with
$\mf x_0=(3,2,-2)^\top$ and $\eta=10^{-10}$.
In Figure~\ref{fig:regularization}(a), we decrease $\alpha$
through $5$, $1$, $0.5$, $10^{-1}$, $3\times10^{-2}$, $10^{-2}$,
$3\times10^{-3}$, $10^{-3}$, $3\times10^{-4}$, $10^{-4}$,
$10^{-5}$, and $10^{-6}$, with $(N,h)=(1024,0.0125)$ and
reference $\alpha_{\rm ref}=10^{-7}$.
In Figure~\ref{fig:regularization}(b), we use
$N\in\{32,64,128,256,512,1024\}$, $h=0.0125$, and
$\alpha\in\{10^{-3},10^{-1},1,5\}$, with $N_{\rm ref}=4096$
for each fixed $\alpha$.
In Figure~\ref{fig:regularization}(c), we halve $h$ from
$0.2$ to $0.0125$, with $N=128$ and the same values of $\alpha$,
using $h_{\rm ref}=0.00078125$ for each fixed $\alpha$.
Panels (b)--(c) report the combined error $E=(E_x^2+E_z^2)^{1/2}$.
The errors decrease under refinement, and larger
$\alpha$ generally yields smaller discretization errors
in these tests, reflecting improved stability
under stronger monotonicity.

\smallskip\noindent{\bf Example 2: a unique solution without strong monotonicity.}
We take $W=W_0$ and retain the distinct recourse constraints.
Since the symmetric part of $\Gamma(\xi)^{-1}$
is positive definite and $R$ has full row rank,
the symmetric part of $M_{22}(\xi)$
is uniformly positive definite.
With the box constraint on $\mf y_1$ imposed directly,
the reduced VI in $(\mf y_1,\mf u)$ is strongly monotone
and has a unique solution.
The first-stage conditions then determine $\lambda_1$
uniquely, giving a unique solution of the two-stage SLCP.

\begin{figure}[htbp]
\centering
\includegraphics[width=0.9\linewidth]{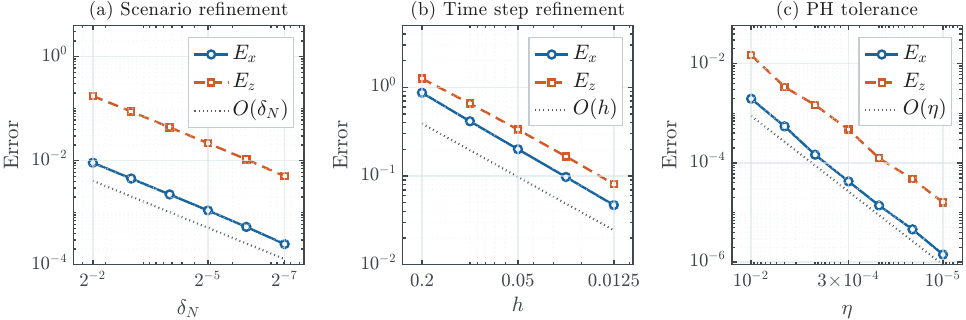}
\caption{Example 2 ($\alpha=0$):
State and decision errors as
(a) $\delta_N$ decreases, (b) $h$ decreases, and
(c) $\eta$ decreases. Dotted lines indicate first-order convergence.}
\label{fig:example2-errors}
\end{figure}

We set $\alpha=0$ and repeat the refinement tests of
Figure~\ref{fig:example1-errors}, keeping all other settings.
Figure~\ref{fig:example2-errors} shows the errors as
$\delta_N$, $h$, and $\eta$ decrease.
The observed orders for $(E_x,E_z)$ are
$(1.0357,1.0217)$, $(1.0501,0.9886)$, and
$(1.0425,0.9780)$, respectively.
Although this example does not satisfy the strong monotonicity
assumption of Theorem~\ref{thm:euler}, the numerical results
show similar convergence behavior.

\smallskip\noindent{\bf Equilibrium and computational results.}
Here we fix $P_t\equiv P=\operatorname{Unif}[-1,1]$ and
$\alpha=0$ in Example 2, so that the system is autonomous.
As in the preceding experiments, we partition $\Xi=[-1,1]$
into $N$ equal subintervals $\Xi_i$.
Since $P$ is uniform, the scenario weights are
$\pi_i=P(\Xi_i)=1/N$, and the discrete distribution is
$
P^N(A)=\frac1N\sum_{i=1}^N\mathbf 1_A(\xi_i),
\  A\in\mathcal B(\Xi).
$

We first examine the convergence to equilibrium using
the left endpoints of $\Xi_i$ as scenario representatives.
We set $N=64$, $h=0.05$, $\eta=10^{-10}$, and $T=10$,
with five initial states:
$\mf x_0^{(1)}=(10,-8,6)^\top$,
$\mf x_0^{(2)}=(-18,15,-16)^\top$,
$\mf x_0^{(3)}=(30,-25,22)^\top$,
$\mf x_0^{(4)}=(-45,38,-35)^\top$, and
$\mf x_0^{(5)}=(65,-55,50)^\top$.

We obtain a common numerical equilibrium reference
$(\mf x^e_{N,\mathrm{ref}},\mf z^e_{N,\mathrm{ref}})$
using the same parameters, with
$\mf x_0=(3,2,-2)^\top$ and $T_{\mathrm{ref}}=35$.
The natural residual is at most $10^{-10}$, and the state
vector field norm is below $10^{-14}$.
We measure the distance to this reference by
\begin{equation}\label{eq:equilibrium-distance}
D(t):=\big(\|x^{\nhe}(t)-\mf x_{N,\mathrm{ref}}^e\|^2
+\|z^{\nhe}(t)-\mf z_{N,\mathrm{ref}}^e\|_{\mc H_{P^N}}^2\big)^{1/2}.
\end{equation}

We next examine how the number of partitions
affects the computed equilibrium payoff.
For each $N\in\{10,50,200,500,1000,3000\}$,
we randomly select one representative $\xi_i$
uniformly from each $\Xi_i$ and repeat the experiment
50 times independently.
The selected representatives remain fixed along each
trajectory, and the coefficient functions remain unchanged.
We set $\mf x_0=(12,-8,10)^\top$, $h=0.1$,
$\eta=10^{-9}$, and a maximum horizon $T_{\max}=30$.
For each run, we compute the natural residual
$\mf r_{\mf x}(\mf z)$ using \eqref{eq:residual-blocks}
with the unregularized operator $\mc M^N$
and the block maximum norm.
At each Euler node, we set
$\tilde F_j^{\nhe}:=N^{-1}\sum_{i=1}^N
\Phi(\mf x_j^{\nhe},\tilde{\mf y}_j^{\nhe},
\tilde{\mf u}_{j,i}^{\nhe},\xi_i)$.
The terminal step $J_N$ is the first $j$ satisfying
$
\max\{\|\tilde F_j^{\nhe}\|,
\|\mf r_{\mf x_j^{\nhe}}(\tilde{\mf z}_j^{\nhe})\|_\infty\}
\le10^{-7}.
$
All runs terminate at $J_N=180$ ($t=18$).
The terminal payoff is defined by
\begin{equation}\label{eq:reported-value}
V_N:=F_1(\mf x^{\nhe}_{J_N},\tilde{\mf y}^{\nhe}_{1,J_N})
+\frac1N\sum_{i=1}^N
F_2(\tilde{\mf x}_{2,J_N,i},
\tilde{\mf y}_{2,J_N,i},\xi_i),
\end{equation}
where the recourse actions
$(\tilde{\mf x}_{2,J_N,i},\tilde{\mf y}_{2,J_N,i})$
are recovered from $\tilde{\mf u}_{J_N,i}^{\nhe}$
using \eqref{eq:primal-recovery}.
Thus, $V_N$ is the payoff at the approximate equilibrium
of the finite-scenario game.

\begin{figure}[htbp]
\centering
\includegraphics[width=0.9\linewidth]{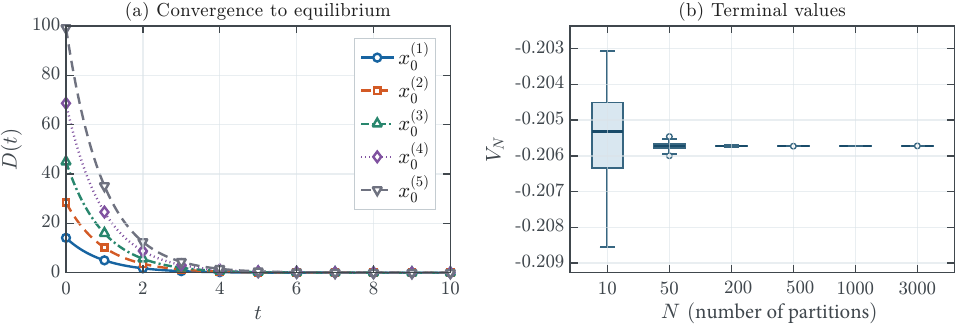}
\caption{Example 2: (a) distance to a common numerical equilibrium from
five initial states; (b) terminal payoffs $V_N$ over 50 runs for each $N$.}
\label{fig:equilibrium-values}
\end{figure}

Figure~\ref{fig:equilibrium-values}(a) shows that all
trajectories approach the same numerical equilibrium.
Figure~\ref{fig:equilibrium-values}(b) shows that
the terminal payoffs become less dispersed as
the number of equal subintervals $N$ increases.
Table~\ref{tab:example2-cost} reports the average results
over 50 runs. The number of PH iterations includes
the final computation. As $N$ increases, the number
of PH iterations remains nearly unchanged, while
the computation time increases.

\begin{table}[htbp]
\centering\small
\caption{Computational results in Example 2 over 50 runs
(mean, with sample variance in parentheses).}
\label{tab:example2-cost}
\setlength{\tabcolsep}{3pt}
\begin{tabular}{@{}rcc@{\qquad}rcc@{}}
\toprule
$N$ & Total PH iter. & Time (s)
& $N$ & Total PH iter. & Time (s)\\
\midrule
10
& $2811.04\ (345.06)$
& $0.101\ (8.66\times10^{-6})$
& 500
& $2857.42\ (1.11)$
& $0.980\ (1.13\times10^{-4})$\\
50
& $2846.04\ (16.73)$
& $0.140\ (2.79\times10^{-5})$
& 1000
& $2858.80\ (0.73)$
& $2.187\ (2.68\times10^{-3})$\\
200
& $2855.26\ (1.26)$
& $0.570\ (6.05\times10^{-5})$
& 3000
& $2859.92\ (0.08)$
& $5.639\ (1.71\times10^{-2})$\\
\bottomrule
\end{tabular}
\end{table}

\section{Conclusion}\label{sec:conclusion}

In this paper, we study the differential two-stage SLCP in~\eqref{DSVI} under a time-dependent
probability distribution. Under strong monotonicity with a scenario-dependent modulus
and suitable joint moment conditions,
we establish Lipschitz continuity in the state and
local $\frac{1}{2}$-H\"older continuity in time for the
solution of the two-stage SLCP, and prove the existence
and uniqueness of a classical solution to the differential two-stage SLCP. We first develop a finite-scenario approximation and
establish an $O(\delta_N)$ convergence rate, where
$\delta_N$ is the mesh size of the scenario partition.
We then apply the forward Euler method with time step $h$
to the finite-scenario problem and finally solve the
resulting finite-dimensional problems by inexact progressive
hedging with stopping tolerance $\eta$.
We establish the total error bound $C(\delta_N+h+\eta)$ of the approximation solution in~\eqref{eq:interpolants} to the  
classical solution of~\eqref{DSVI} in Theorem~\ref{thm:euler}, where
$C$ is independent of $N$, $h$ and $\eta$.
Finally, we apply the method to a quadratic two-stage
stochastic two-player zero-sum game and relate the
equilibria of the resulting differential two-stage SLCP to the KKT points
of its minimax formulation. Numerical results support the first-order error
bounds under regularization and show similar convergence
behavior without strong monotonicity. Future research may extend the analysis to monotone two-stage
SLCPs without strong monotonicity and investigate the
long-time convergence of trajectories to equilibria.

\vspace{-0.1in}
\section*{Acknowledgments}

The authors used AI tools to assist in searching for
related references and improving the readability of
this manuscript. All AI-assisted content was carefully
reviewed and revised by the authors.
The authors assume responsibility for all content.

\bibliographystyle{amsplain}
\bibliography{references_SICON}

\end{document}